\documentclass[11pt]{article}
\usepackage[letterpaper,margin=1in]{geometry}
\usepackage{amsmath,amssymb,amsthm,mathtools,mathrsfs}
\usepackage[T1]{fontenc}
\usepackage{lmodern,microtype}
\usepackage{xcolor}
\usepackage{tikz,pgfplots}
\pgfplotsset{compat=1.3}
\usepackage{cite}
\usepackage[colorlinks=true,linkcolor=blue,citecolor=blue,urlcolor=blue]{hyperref}
\hypersetup{pdftitle={Threshold Attainment and Sharp Asymptotics for Lp-Perturbed Moser-Trudinger Functionals on R2},
 pdfsubject={Threshold extremals, strong-penalty asymptotics and profile selection},
 pdfauthor={Kaiwen Guo and Yanjun Liu}}
\usepackage{enumitem}
\usepackage{needspace}
\usepackage{placeins}
\setlist[enumerate]{label=\textup{(\roman*)},itemsep=3pt,topsep=5pt}
\DeclareRobustCommand{\reftext}[2]{#1~\ref{#2}}
\renewcommand{\eqref}[1]{\textup{(\ref{#1})}}
\makeatletter
\renewenvironment{proof}[1][\proofname]{\par
  \pushQED{\qed}\normalfont\topsep6\p@\@plus6\p@\relax
  \trivlist\item[\hskip\labelsep\bfseries #1\@addpunct{.}]\ignorespaces
}{\popQED\endtrivlist\@endpefalse}
\makeatother
\allowdisplaybreaks[2]
\newtheorem{theorem}{Theorem}[section]
\newtheorem{lemma}[theorem]{Lemma}
\newtheorem{proposition}[theorem]{Proposition}
\newtheorem{corollary}[theorem]{Corollary}
\theoremstyle{remark}
\newtheorem{remark}[theorem]{Remark}
\newcommand{\R}{\mathbb R}
\newcommand{\A}{\mathcal A}
\newcommand{\Cp}{\mathcal C_p}
\newcommand{\Hfun}{\mathcal H_p}
\newcommand{\lamstar}{\lambda_*(p)}
\newcommand{\dd}{\,dx}
\newcommand{\weak}{\rightharpoonup}
\numberwithin{equation}{section}
\title{Threshold Attainment and Sharp Asymptotics
for $L^p$-Perturbed Moser--Trudinger Functionals on $\R^2$}
\author{Kaiwen Guo$^{a}$ \qquad Yanjun Liu$^{b}$\thanks{Corresponding author.}\\[5pt]
\small $^{a}$Mathematical Research Center, Shandong University,\\
\small Jinan 250100, P. R. China\\
\small \href{mailto:kwguo@mail.nankai.edu.cn}{\textcolor{black}{kwguo@mail.nankai.edu.cn}}\\[4pt]
\small $^{b}$School of Mathematical Sciences, Chongqing Normal University,\\
\small Chongqing 401331, P. R. China\\
\small \href{mailto:liuyj@mail.nankai.edu.cn}{\textcolor{black}{liuyj@mail.nankai.edu.cn}}}
\date{}

\begin{document}
\maketitle
\begin{abstract}
For the critical Moser--Trudinger supremum
\[
 S(\lambda,p)=\sup_{\|u\|_{H^1(\R^2)}\le1}
 \int_{\R^2}\bigl(e^{4\pi u^2}-1-\lambda|u|^p\bigr)\,dx,
 \qquad p\ge2,
\]
Chen--Jiang--Lu--Zhu gave the first classification of attainment away
from the thresholds across the full range of $p$. We prove attainment
at both thresholds for $p=2$ and at the unique threshold for each
$2<p\le4$. At and above the upper quadratic threshold,
$S(\lambda,2)=4\pi e^{1-2\gamma_{\mathrm E}}$, Nguyen's concentration
level, where $\gamma_{\mathrm E}$ is the Euler--Mascheroni constant.
At each endpoint, the centred radial extremals form a compact set and
their $L^p$ masses determine the slope from within the attainment interval.
For fixed $p>4$, we obtain the expansion
\[
 S(\lambda,p)=4\pi+\mathcal C_p\lambda^{-2/(p-4)}
 +D_p\lambda^{-4/(p-4)}+o\bigl(\lambda^{-4/(p-4)}\bigr),
 \qquad \lambda\to\infty.
\]
Here $\mathcal C_p$ is the maximum of a limiting problem at unit $L^2$
mass, and $D_p$ selects the limiting profiles. For fixed $\lambda$
above the quartic threshold $\lambda_*(4)$, we prove
\[
 \lim_{p\downarrow4}(p-4)\log\bigl(S(\lambda,p)-4\pi\bigr)
 =2\log\frac{\lambda_{\mathrm{loc}}}{\lambda}<0,
\]
where $\lambda_{\mathrm{loc}}$ is determined by the sharp
Gagliardo--Nirenberg constant. An exact radial transformation extends
these results to every planar norm.
\end{abstract}



\section{Introduction and main results}
We examine how an additive $L^p$ penalty changes attainment and
the profiles of extremals for the critical Moser--Trudinger supremum
\begin{equation}
 \begin{aligned}
 \A&=\left\{u\in H^1(\R^2):
       \int_{\R^2}(|\nabla u|^2+u^2)\,dx\le1\right\},\\
 J_{\lambda,p}(u)&=\int_{\R^2}
       \bigl(e^{4\pi u^2}-1-\lambda|u|^p\bigr)\,dx,\\
 S(\lambda,p)&=\sup_{u\in\A}J_{\lambda,p}(u),\qquad p\ge2.
 \end{aligned}
 \label{eq:problem}
\end{equation}
The quadratic penalty changes the vanishing barrier itself. For $p>2$,
the competition between the quartic exponential term and the $p$th-power
penalty singles out $p=4$ as the transition exponent.
After radial rearrangement, compactness can
fail in two ways: concentration forms an increasingly narrow peak,
whereas vanishing spreads the $L^2$ mass over larger spatial scales.
We determine whether extremals survive at the transition parameters
and how they spread out when the penalty becomes large or the exponent
approaches $4$ from above.
All  norms are taken over $\R^2$, and $B_R$ denotes the ball
centred at the origin with radius $R$. We use $\|u\|_{H^1}^2=\|\nabla u\|_2^2+\|u\|_2^2$
and call a function normalized when $\|u\|_{H^1}=1$. By a centred
radial function we mean a nonnegative function that is radial about the
origin and nonincreasing in the radius; the compactness statements below
are formulated for these representatives.

The critical exponential embedding was established by Trudinger and
Moser~\cite{Trudinger,Moser}. Its extremal problem already exhibits a
competition between concentration and compactness on bounded domains.
Carleson--Chang~\cite{CarlesonChang} proved attainment on the disk,
and Flucher~\cite{Flucher} extended this result to general planar domains
by Green-function transplantation. Lions' concentration-compactness
principle~\cite{Lions} provides the framework for this analysis, while
de Figueiredo--do \'O--Ruf~\cite{deFigueiredoDoORuf} developed
concentration tests for more general critical nonlinearities. The refined
energy expansions of Mancini--Martinazzi~\cite{ManciniMartinazzi} and
the perturbation criteria of Thizy~\cite{Thizy} show how lower-order
terms can decide whether the critical value is attained.

On unbounded domains, subcritical exponential estimates were obtained
by Cao~\cite{Cao}, do \'O~\cite{doO1997} and
Adachi--Tanaka~\cite{AdachiTanaka}. Ruf~\cite{Ruf} proved the sharp
planar inequality under the full Sobolev norm, and Li--Ruf~\cite{LiRuf}
established the higher-dimensional counterpart and attainment.
Relations with the subcritical theory were developed by
Cassani--Sani--Tarsi~\cite{CassaniSaniTarsi} and
Lam--Lu--Zhang~\cite{LamLuZhang2017,LamLuZhang2019}.
Unlike the bounded Dirichlet setting, the full norm permits vanishing
at unit mass. Its role in nonattainment was identified by
Ishiwata~\cite{Ishiwata}, and a vanishing-concentration-compactness
alternative was proved by do \'O--Sani--Tarsi~\cite{doOSaniTarsi}.
Ikoma--Ishiwata--Wadade~\cite{IkomaIshiwataWadade} studied attainment
under inhomogeneous norm constraints, and Nguyen~\cite{Nguyen}
identified the sharp normalized vanishing and concentration levels.
Under a homogeneous gradient constraint, the exact-growth results of
Ibrahim--Masmoudi--Nakanishi~\cite{IMN} and the threshold criteria of
Ibrahim--Masmoudi--Nakanishi--Sani~\cite{IMNS} likewise show that both
the behavior near zero and the admissible norm govern compactness.

For additive $L^p$ perturbations, Hashizume~\cite{Hashizume2020}
obtained a threshold classification on the disk and subsequently treated
more general lower-order terms~\cite{HashizumeLowerOrder}.
Chen--Jiang--Lu--Zhu~\cite{ChenBounded} established the corresponding
classification away from the threshold on bounded planar domains.
On the whole plane, Chen--Jiang--Lu--Zhu~\cite{Chen} gave, to our
knowledge, the first systematic classification of attainment for additive
$L^p$ perturbations in \eqref{eq:problem} across the full range $p\ge2$.
For $p=2$, Chen--Lu--Zhu~\cite[Theorem 1.7]{Chen2020} had proved
attainment for $\lambda_-<\lambda<4\pi$.
Chen--Jiang--Lu--Zhu~\cite[Theorem 1.2]{Chen} proved attainment for
$4\pi\le\lambda<\lambda_+$. Together these results give an extremal for
every $\lambda\in(\lambda_-,\lambda_+)$, as stated in their
Corollary~1.3. For $2<p\le4$, their Theorem~1.1 gives attainment for
every $\lambda<\lambda_*(p)$ and nonattainment for every
$\lambda>\lambda_*(p)$. For $p>4$, the same theorem gives an extremal
for every $\lambda\in\R$. They also identified concentration above the
upper quadratic threshold and vanishing above the threshold for $2<p\le4$
as the mechanisms for nonattainment.
We settle the endpoint cases left open by this classification by proving attainment at
$\lambda=\lambda_-$ and $\lambda=\lambda_+$ for $p=2$, and at
$\lambda=\lambda_*(p)$ for $2<p\le4$. We also prove compactness of centred radial extremals
at these endpoints and identify the extremal $L^p$ mass that determines
each one-sided slope.

The two quadratic endpoints require different estimates. At the lower
endpoint, the Gagliardo--Nirenberg bound appearing in
Chen--Lu--Zhu~\cite[Theorem 1.7]{Chen2020} must be made strict to exclude
vanishing. We use an algebraic test together with sharp interpolation
estimates, whose optimizer formulation and Fourier bounds go back to
Weinstein~\cite{Weinstein} and Beckner~\cite{Beckner}. At the upper
endpoint, concentration must be excluded at the sharp normalized level
\begin{equation}
 C_{\mathrm{cl}}=4\pi e^{1-2\gamma_{\mathrm E}}<4\pi,
 \label{eq:concentration}
\end{equation}
where $\gamma_{\mathrm E}=\lim_{n\to\infty}
(\sum_{k=1}^n k^{-1}-\log n)$ is the Euler--Mascheroni constant.
Chen--Jiang--Lu--Zhu~\cite[Theorem 1.2]{Chen} give $\pi e$ as the value
of the quadratic supremum above the upper threshold. Immediately after
equation~(3.7) in version~1 of their preprint, they evaluate $\pi e^{1+4\pi A}$ as
$\pi e$, where $A$ is the regular part at the origin of the Green
function for $-\Delta+1$. As Nguyen observed after Theorem~1.5
of~\cite{Nguyen}, the normalization $(-\Delta+1)G=\delta_0$ gives
$A=(\log 2-\gamma_{\mathrm E})/(2\pi)$, so
$\pi e^{1+4\pi A}=4\pi e^{1-2\gamma_{\mathrm E}}>\pi e$.
The Green expansion is recorded in \eqref{eq:green-normalization}.
For $p>2$, the comparison used in the proof of
\cite[Theorem 1.1]{Chen} still excludes concentration, since
$S(\lambda,p)\ge4\pi>C_{\mathrm{cl}}$.
Since the $L^2$ mass of a concentrating sequence tends to zero, a fixed
quadratic penalty preserves its limiting value. We prove that
$C_{\mathrm{cl}}$ is the value of $S(\lambda,2)$ at and above the
upper threshold. A fourth-order expansion excludes concentration of
extremals at this endpoint. At $p=4$,
a sixth-moment estimate for a sharp Gagliardo--Nirenberg optimizer
provides the strict gap needed to rule out vanishing.

Our second contribution concerns $p>4$, where attainment at every finite
parameter is already known from~\cite{Chen}. We determine the rate at
which $S(\lambda,p)$ approaches $4\pi$ as $\lambda\to\infty$, including
the next term, and describe the profiles of its spreading extremals.
After rescaling, these profiles converge strongly along subsequences to
maximizers of a limiting functional at unit $L^2$ mass. The second
coefficient combines the sixth moment with the correction required by
the full $H^1$ constraint and selects the possible limiting maximizers.
Uniform remainder estimates along the varying rescaled profiles justify
the two-term expansion.
Related fixed-mass problems occur in the cubic-quintic model of
Li--Wang~\cite{LiWang}, and vanishing under parameter-dependent scaling
is studied by Hashizume~\cite{Hashizume2025}.
Hashizume~\cite{Hashizume2023} also treats a subcritical Neumann
functional with parameter-dependent normalization on a bounded domain.
Here the norm constraint is fixed and the parameter enters only the
objective, which makes the normalization correction part of the
second asymptotic coefficient.

We further quantify the transition from attainment for $p>4$ to
nonattainment at $p=4$. For a fixed $\lambda>\lambda_*(4)$, the excess
value above $4\pi$ and the gradient energy decay exponentially in
$1/(p-4)$, with a rate determined by the sharp Gagliardo--Nirenberg
constant. After normalization of both mass and gradient energy, the
extremals converge along subsequences to its optimizers. This result
describes how the change of regime in~\cite{Chen} occurs and complements
the expansion at a fixed exponent. Its proof uses estimates uniform
in $p$ and does not exchange the limits $\lambda\to\infty$ and
$p\downarrow4$.

The anisotropic extension is based on convex symmetrization
\cite{AlvinoEtAl} and the Finsler--Laplacian framework
\cite{FeroneKawohl}. Wang--Xia~\cite{WangXia} developed the planar
Finsler--Liouville blow-up theory, and Zhou--Zhou proved
Moser--Trudinger inequalities and attainment on bounded domains
\cite{ZhouZhou2018,ZhouZhou2019} and under the full anisotropic Sobolev
norm on unbounded domains~\cite{ZhouZhou2020}. Singular weighted
whole-space inequalities were established by
Adimurthi--Yang~\cite{AdimurthiYang} and, for anisotropic norms,
Guo--Liu~\cite{GuoLiuEntire}. Liu developed singular anisotropic
concentration-compactness~\cite{Liu2020Concentration} and studied
exact-growth inequalities and maximizers~\cite{Liu2022Anisotropic}.
Guo--Liu~\cite{GuoLiuWulff} obtained existence and Wulff symmetry,
and their recent work~\cite{GuoLiuSingularThresholds} addresses critical
attainment and thresholds under singular anisotropic mixed-norm
constraints. On bounded domains, the $L^p$ improvement of
Liu~\cite{LiuImproved} and its singular counterpart due to
Guo--Liu~\cite{GuoLiuSingularBounded} modify the admissible norm.
The extension proved here instead retains a fixed anisotropic norm and
an additive penalty. An exact radial transfer carries the endpoint and
asymptotic conclusions to every constant planar norm, with a direct
proof in Section~\ref{sec:strong-proofs}.

% Original reference 20 and its discussion are temporarily omitted.
% Our companion manuscript~\cite{GuoLiu} studies bounded planar domains
% under a Dirichlet constraint with smooth, strongly convex anisotropy.
% It proves threshold attainment and power-mass selection for
% $1\le p\le2$ and describes concentration for $p>2$, with corrections
% and concentration points determined by Green moments; see
% \cite[Theorems 1.2--1.4]{GuoLiu}. The full $H^1$ constraint in the
% present paper adds the vanishing barrier and moves the transition to
% $p=4$: for $p>4$, profiles spread out and converge, after rescaling, to
% maximizers of a limiting problem at unit $L^2$ mass. We also give a short proof of the whole-space transfer identity
% from \cite[Theorem 1.5(ii)]{GuoLiu}, so that the anisotropic conclusions
% follow explicitly from the Euclidean results.

For the quadratic penalty, define
\begin{equation}
 \begin{aligned}
 D_*&:=\sup_{u\in\A\setminus\{0\}}
 \frac{\int_{\R^2}(e^{4\pi u^2}-1-4\pi u^2)\dd}{\|\nabla u\|_2^2},
 &\lambda_-&:=4\pi-D_*,\\
 \lambda_+&:=\sup_{u\in\A\setminus\{0\}}
 \frac{\int_{\R^2}(e^{4\pi u^2}-1)\dd-C_{\mathrm{cl}}}{\|u\|_2^2}.
 \end{aligned}
 \label{eq:quadratic-quotients}
\end{equation}

Chen--Jiang--Lu--Zhu~\cite[Corollary 1.3]{Chen} give the open
interval of attainment. The following result determines the values
and the extremals at its endpoints.

\begin{theorem}\label{thm:quadratic}
The quadratic thresholds satisfy
\begin{equation}
 4\pi<D_*<\infty,\qquad
 -\infty<\lambda_-<0,\qquad
 4\pi+C_{\mathrm{cl}}<\lambda_+<\infty.
 \label{eq:quadratic-threshold-bounds}
\end{equation}
\begin{enumerate}
\item
Both quotients in \eqref{eq:quadratic-quotients} attain their maxima.
The supremum is attained at $\lambda_-$ and $\lambda_+$, with
\[
 S(\lambda_-,2)=D_*,\qquad
 S(\lambda,2)=C_{\mathrm{cl}}\quad\text{for all }\lambda\ge\lambda_+.
\]

\item At either endpoint, centred radial extremals form a compact set
in $H^1$. Those approaching either endpoint from inside the interval
are precompact in $H^1$, and every strong limit is an endpoint
extremal. Normalized centred radial maximizing sequences for the
first quotient in \eqref{eq:quadratic-quotients} are precompact in $H^1$.

The endpoint masses give the one-sided slopes. Let $m_-$ be the
minimum of $\|u\|_2^2$ among extremals at
$\lambda_-$, and $m_+$ the maximum among extremals at $\lambda_+$.
Then $0<m_-<1$, $0<m_+<1$, and
\begin{align}
 S(\lambda_-+h,2)&=D_*-m_-h+o(h),\label{eq:quadratic-lower-slope}\\
 S(\lambda_+-h,2)&=C_{\mathrm{cl}}+m_+h+o(h),\qquad h\downarrow0.
 \label{eq:quadratic-upper-slope}
\end{align}
Every family of extremals at these parameters has $L^2$ masses
converging to $m_-$ and $m_+$, respectively; every strong limit of
these centred radial extremals has the corresponding mass. Hence
\[
 S'_-(\lambda_-,2)=-1,\quad S'_+(\lambda_-,2)=-m_-,\quad
 S'_-(\lambda_+,2)=-m_+,\quad S'_+(\lambda_+,2)=0.
\]
\end{enumerate}
\end{theorem}

For $2<p\le4$, let $\lamstar\in(0,\infty)$ denote the threshold
of Chen--Jiang--Lu--Zhu~\cite[Theorem 1.1]{Chen}.

\Needspace{12\baselineskip}
\begin{theorem}\label{thm:endpoint}
Let $2<p\le4$.
\begin{enumerate}
\item The supremum $S(\lamstar,p)=4\pi$ is attained, and
\begin{equation}
 \lamstar=\max_{u\in\A\setminus\{0\}}
 \frac{\int_{\R^2}(e^{4\pi u^2}-1)\dd-4\pi}{\|u\|_p^p}.
 \label{eq:threshold-quotient}
\end{equation}

\item Centred radial endpoint extremals form a compact set in $H^1$.
If $\lambda_n<\lamstar$, $\lambda_n\to\lamstar$, and $u_n$ are centred
radial extremals for $S(\lambda_n,p)$, a subsequence converges strongly
in $H^1$ to an endpoint extremal. The maximal endpoint mass determines
the one-sided slope. Set
\begin{equation}
 \sigma_p=\max\{\|u\|_p^p:u\in\A,\ J_{\lambda_*(p),p}(u)=4\pi\}>0.
 \label{eq:sigma}
\end{equation}
For every family of extremals $u_h$ at $\lamstar-h$,
\begin{align}
 S(\lamstar-h,p)&=4\pi+\sigma_p h+o(h),
 \label{eq:endpoint-value}\\
 \|u_h\|_p^p&\longrightarrow\sigma_p,\qquad h\downarrow0.
 \label{eq:endpoint-mass}
\end{align}
Thus $S'_-(\lamstar,p)=-\sigma_p$. Every strong limit of centred radial
$u_h$ as $h\downarrow0$ is an endpoint extremal with power mass $\sigma_p$.
\end{enumerate}
\end{theorem}

For $p>4$, Chen--Jiang--Lu--Zhu~\cite[Theorem 1.1]{Chen}
proved that extremals exist for every finite $\lambda$.
We determine how these extremals spread as $\lambda\to+\infty$. The leading
correction to $4\pi$ is governed by the unit-mass problem
\begin{equation}
 \Hfun(v)=8\pi^2\|v\|_4^4-4\pi\|\nabla v\|_2^2-\|v\|_p^p,\qquad
 \Cp=\sup_{\substack{v\in H^1(\R^2)\\\|v\|_2=1}}\Hfun(v).
 \label{eq:limit-problem}
\end{equation}

Let $\mathcal G_p$ be the set of centred radial maximizers. The second
correction includes the sixth moment and the full-norm normalization:
\begin{equation}
 D_p=\max_{v\in\mathcal G_p}
 \left\{\frac{32\pi^3}{3}\|v\|_6^6-4\pi \|\nabla v\|_2^4-\Cp \|\nabla v\|_2^2\right\}.
 \label{eq:second-coefficient}
\end{equation}
The compactness in Theorem~\ref{thm:strong} makes this maximum
finite and attained.

\Needspace{10\baselineskip}
\begin{theorem}\label{thm:strong}
Fix $p>4$.
\begin{enumerate}
\item The value $\Cp$ in \eqref{eq:limit-problem} is finite and
positive, and
\begin{equation}
S(\lambda,p)=4\pi+\Cp\lambda^{-2/(p-4)}
 +D_p\lambda^{-4/(p-4)}
 +o\bigl(\lambda^{-4/(p-4)}\bigr),\qquad\lambda\to+\infty,
 \label{eq:main-asymptotic}
\end{equation}
where $D_p$ is given by \eqref{eq:second-coefficient}; its sign is not
prescribed. The corresponding moments satisfy, for every family of
extremals $u_\lambda$,
\begin{align}
 \lambda^{(p-2)/(p-4)}\|u_\lambda\|_p^p
 &\longrightarrow\frac{2\Cp}{p-4},\label{eq:strong-mass}\\
\lambda^{2/(p-4)}\bigl[\int_{\R^2}(e^{4\pi u_\lambda^2}-1)\dd-4\pi\bigr]
 &\longrightarrow\frac{p-2}{p-4}\Cp,\qquad\lambda\to+\infty.
 \label{eq:strong-exp}
\end{align}

\item The unit-mass problem \eqref{eq:limit-problem} has a centred
radial maximizer, and every centred radial maximizing sequence is
precompact in $H^1$. The second coefficient also selects the rescaled
profiles. For any family of centred radial extremals
$u_\lambda$, set
\begin{equation}
 v_\lambda(x)=\lambda^{1/(p-4)}
 u_\lambda\bigl(\lambda^{1/(p-4)}x\bigr),\qquad\lambda>0.
 \label{eq:profile-scaling}
\end{equation}
Every sequence $\lambda_n\to+\infty$ has a subsequence for which
$v_{\lambda_n}\to v$ strongly in $H^1$, where $v\in\mathcal G_p$ and
\begin{equation}
 \frac{32\pi^3}{3}\|v\|_6^6-4\pi \|\nabla v\|_2^4-\Cp \|\nabla v\|_2^2=D_p.
 \label{eq:second-profile-selection}
\end{equation}
\end{enumerate}
\end{theorem}

The expansion is for each fixed $p>4$ and requires neither
uniqueness nor nondegeneracy of the limiting maximizers.

\begin{figure}[!htbp]
\centering
\begin{minipage}[t]{0.485\textwidth}
\centering
\begin{tikzpicture}[x=1.03cm,y=0.86cm,font=\footnotesize,
                    line cap=round,line join=round]
 \definecolor{regimeblue}{RGB}{35,83,125}
 \draw[->,black] (0.35,0.6)--(6.75,0.6) node[right] {$\lambda$};
 \draw[->,black] (0.35,0.6)--(0.35,5.85);
 \node[anchor=south west] at (0.40,5.47) {$S(\lambda,2)$};
 \node[font=\footnotesize\bfseries] at (3.50,6.18) {(a) $p=2$};
 % The affine barrier has a steeper slope than the attained branch.
 \draw[gray!65,densely dashed] (2,3.85)--(4.95,0.605);
 \draw[gray!65,densely dashed] (0.35,1.60)--(5,1.60);
 \node[anchor=east] at (0.28,1.60) {$C_{\mathrm{cl}}$};
 \node[anchor=south] at (1.56,4.62) {$4\pi-\lambda$};
 \draw[gray!80,line width=1.1pt] (0.65,5.335)--(2,3.85);
 % Bezier control slopes -1.00,-0.80,-0.45 ensure convexity.
 \draw[regimeblue,line width=1.3pt]
       (2,3.85)..controls (3,2.85) and (4,2.05)..(5,1.60);
 \draw[gray!80,line width=1.1pt] (5,1.60)--(6.55,1.60);
 \draw[gray!65,densely dotted] (2,0.6)--(2,3.85);
 \draw[gray!65,densely dotted] (5,0.6)--(5,1.60);
 \fill[regimeblue] (2,3.85) circle (2.3pt);
 \fill[regimeblue] (5,1.60) circle (2.3pt);
 \node[anchor=north] at (2,0.46) {$\lambda_-$};
 \node[anchor=north] at (5,0.46) {$\lambda_+$};
 \draw (3.20,0.64)--(3.20,0.56);
 \node[anchor=north] at (3.20,0.46) {$0$};
 \draw[regimeblue,line width=0.8pt] (2,-0.02)--(5,-0.02);
 \draw[regimeblue,line width=0.8pt] (2,0.06)--(2,-0.10);
 \draw[regimeblue,line width=0.8pt] (5,0.06)--(5,-0.10);
 \node[regimeblue,anchor=north] at (3.5,-0.20) {attained on $[\lambda_-,\lambda_+]$};
\end{tikzpicture}
\end{minipage}\hfill
\begin{minipage}[t]{0.485\textwidth}
\centering
\begin{tikzpicture}[x=1.03cm,y=0.86cm,font=\footnotesize,
                    line cap=round,line join=round]
 \definecolor{regimeblue}{RGB}{35,83,125}
 \definecolor{regimeorange}{RGB}{163,92,35}
 \draw[->,black] (0.35,0.6)--(6.75,0.6) node[right] {$\lambda$};
 \draw[->,black] (0.35,0.6)--(0.35,5.85);
 \node[anchor=south west] at (0.40,5.47) {$S(\lambda,p)$};
 \node[font=\footnotesize\bfseries] at (3.50,6.18) {(b) $p>2$};
 \draw[gray!65,densely dashed] (0.35,1.60)--(4,1.60);
 \node[anchor=east] at (0.28,1.60) {$4\pi$};
 \draw[regimeblue,line width=1.3pt]
       (0.35,4.80)..controls (1.80,3.00) and (3.00,2.20)..(4,1.60);
 \draw[gray!80,line width=1.1pt] (4,1.60)--(6.55,1.60);
 \draw[gray!65,densely dotted] (4,0.6)--(4,1.60);
 \fill[regimeblue] (4,1.60) circle (2.3pt);
 \node[anchor=north] at (4,0.46) {$\lamstar$};
 \node[anchor=north] at (0.35,0.46) {$0$};
 \node[regimeblue,fill=white,inner sep=2pt] at (2.20,1.95) {$2<p\le4$};
 \draw[regimeorange,line width=1.3pt,dash pattern=on 5pt off 2.5pt,->]
       (0.35,4.80)..controls (2.10,3.40) and (3.55,2.40)..(6.55,1.90);
 \node[regimeorange,fill=white,inner sep=2pt] at (4.83,2.73) {$p>4$};
 \node[regimeblue,anchor=north,align=center] at (3.5,-0.20)
       {threshold attained for $2<p\le4$};
\end{tikzpicture}
\end{minipage}
\caption{Attainment regimes (schematic; axes not to scale). Filled endpoints are attained; grey branches are not. For $p>4$ (dashed orange), the value remains above $4\pi$ at every finite parameter and the supremum is attained.}
\label{fig:attainment-regimes}
\end{figure}
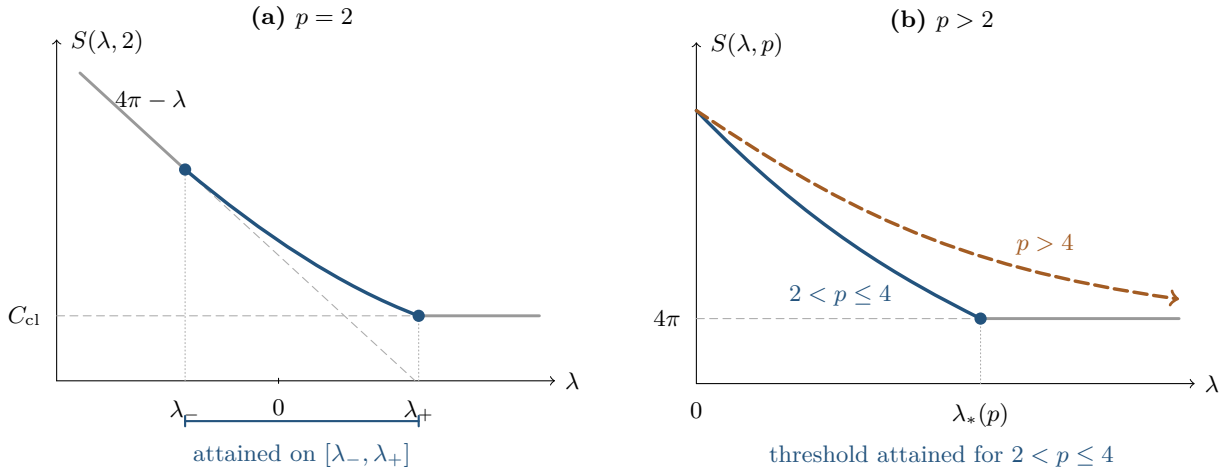

\FloatBarrier
Figure~\ref{fig:attainment-regimes} combines the classification away
from the thresholds in~\cite{Chen} with the attained endpoints of
Theorems~\ref{thm:quadratic} and~\ref{thm:endpoint}.

\FloatBarrier

We next fix the penalty and let $p\downarrow4$. Set
\begin{equation}
 B_2=\sup_{u\in H^1(\R^2)\setminus\{0\}}
       \frac{\|u\|_4^4}{\|\nabla u\|_2^2\|u\|_2^2},\qquad
 \lambda_{\mathrm{loc}}=8\pi^2-\frac{4\pi}{B_2}.
 \label{eq:gn-sharp}
\end{equation}
At $\lambda_{\mathrm{loc}}$, the leading small-gradient gain vanishes;
$\lambda_*(4)$ is the global attainment threshold.

\begin{theorem}\label{thm:p4-transition}
The quartic threshold satisfies
$\lambda_*(4)>\lambda_{\mathrm{loc}}>0$.
Fix $\lambda>\lambda_*(4)$ and write $p=4+\varepsilon$.
\begin{enumerate}
\item As $\varepsilon\downarrow0$, the excess above $4\pi$ satisfies
\begin{equation}
 \lim_{\varepsilon\downarrow0}\varepsilon
 \log\bigl(S(\lambda,4+\varepsilon)-4\pi\bigr)
 =2\log\frac{\lambda_{\mathrm{loc}}}{\lambda}.
\label{eq:p4-value-rate}
\end{equation}
For every family of extremals $u_\varepsilon$, put
$g_\varepsilon=\|\nabla u_\varepsilon\|_2^2$. Then
$g_\varepsilon\to0$, $\|u_\varepsilon\|_2^2=1-g_\varepsilon\to1$, and
\begin{equation}
 \begin{aligned}
 \lim_{\varepsilon\downarrow0}\varepsilon\log g_\varepsilon
 &=2\log\frac{\lambda_{\mathrm{loc}}}{\lambda},\\
 \lim_{\varepsilon\downarrow0}\varepsilon
 \log\sqrt{\frac{1-g_\varepsilon}{g_\varepsilon}}
 &=\log\frac{\lambda}{\lambda_{\mathrm{loc}}}.
 \end{aligned}
\label{eq:p4-scale-rate}
\end{equation}

\item To identify the profiles, normalize mass and gradient energy.
For centred radial extremals, set
\begin{equation}
 v_\varepsilon(y)=g_\varepsilon^{-1/2}
 u_\varepsilon\!\left(
 \sqrt{\frac{1-g_\varepsilon}{g_\varepsilon}}\,y\right).
 \label{eq:p4-profile-scaling}
\end{equation}
Every sequence $\varepsilon_n\downarrow0$ has a subsequence for which
$v_{\varepsilon_n}\to V$ strongly in $H^1(\R^2)$, where $V$ is a
centred radial Gagliardo--Nirenberg optimizer satisfying
\[
 \|V\|_2=\|\nabla V\|_2=1,\qquad \|V\|_4^4=B_2.
\]
\end{enumerate}
\end{theorem}

Thus the spreading scale grows exponentially as $p\downarrow4$.
This rate is determined by $\lambda_{\mathrm{loc}}$, which is
strictly below the attainment threshold $\lambda_*(4)$.

Let $F$ be any norm on $\R^2$, let $F^\circ$ be its polar norm, and put
\[
 \kappa_2=|\{x\in\R^2:F^\circ(x)<1\}|,\qquad
 \theta=\frac{\kappa_2}{\pi}.
\]
Here $\kappa_2$ is the area of the unit Wulff ball. For the Euclidean
norm, $\kappa_2=\pi$ and $\theta=1$.
A centred Wulff-radial function has the form $u(x)=U(F^\circ(x))$,
where $U$ is nonnegative and nonincreasing.
Define
\begin{equation}
 S_F(\lambda,p)=
 \sup_{\substack{u\in H^1(\R^2)\\\int_{\R^2}(F(\nabla u)^2+|u|^2)\dd\le1}}
 \int_{\R^2}\bigl(e^{4\kappa_2 u^2}-1-\lambda|u|^p\bigr)\dd.
 \label{eq:anisotropic-problem}
\end{equation}

The exact relation is
\begin{equation}
 S_F(\lambda,p)=\theta S(\lambda\theta^{-p/2},p).
 \label{eq:transfer}
\end{equation}
The anisotropic supremum is attained if and only if the Euclidean
supremum at parameter $\lambda\theta^{-p/2}$ is attained.

\begin{corollary}\label{cor:anisotropy}
For $2<p\le4$, the anisotropic threshold is attained and
\begin{equation}
 \lambda_*^F(p)=\theta^{p/2}\lamstar.
 \label{eq:anisotropic-endpoint}
\end{equation}
As $h\downarrow0$,
\[
 S_F(\lambda_*^F(p)-h,p)
 =4\kappa_2+\theta^{1-p/2}\sigma_p h+o(h).
\]
For $p>4$, as $\lambda\to\infty$, one has
\begin{equation}
 S_F(\lambda,p)=4\kappa_2+
 \theta^{2(p-2)/(p-4)}\Cp\lambda^{-2/(p-4)}
 +\theta^{(3p-4)/(p-4)}D_p\lambda^{-4/(p-4)}
 +o\bigl(\lambda^{-4/(p-4)}\bigr).
 \label{eq:anisotropic-asymptotic}
\end{equation}
For any family of anisotropic extremals,
\begin{align*}
 \lambda^{(p-2)/(p-4)}\|u_\lambda\|_p^p&\longrightarrow
 \frac{2\theta^{2(p-2)/(p-4)}\Cp}{p-4},\\
 \lambda^{2/(p-4)}
 \left[\int_{\R^2}(e^{4\kappa_2 u_\lambda^2}-1)\dd-4\kappa_2\right]
 &\longrightarrow\frac{p-2}{p-4}\theta^{2(p-2)/(p-4)}\Cp.
\end{align*}
For fixed $\lambda>\theta^2\lambda_*(4)$, the critical transition is
\begin{equation}
 \lim_{p\downarrow4}(p-4)
 \log\bigl(S_F(\lambda,p)-4\kappa_2\bigr)
 =2\log\frac{\theta^2\lambda_{\mathrm{loc}}}{\lambda}.
 \label{eq:anisotropic-p4-rate}
\end{equation}
\end{corollary}

\begin{corollary}\label{cor:quadratic-transfer}
For the norm $F$ above, the quadratic barriers and thresholds satisfy
\[
 C_{\mathrm{cl},F}=\theta C_{\mathrm{cl}},\qquad
 \lambda_-^F=\theta\lambda_-,\qquad
 \lambda_+^F=\theta\lambda_+.
\]
In particular, $\lambda_+^F>\theta(4\pi+C_{\mathrm{cl}})$ is finite and
$S_F(\theta(4\pi+C_{\mathrm{cl}}),2)>C_{\mathrm{cl},F}$.
Both endpoints are attained, with compact centred Wulff-radial
extremal sets. Let $m_-$ and $m_+$ denote, respectively, the minimal
and maximal $L^2$ masses of Euclidean extremals at $\lambda_-$ and
$\lambda_+$, as in Theorem~\ref{thm:quadratic}.
One has $\lambda_-^F<0$ and, as $h\downarrow0$,
\begin{align*}
 S_F(\lambda_-^F+h,2)&=\theta D_*-m_-h+o(h),\\
 S_F(\lambda_+^F-h,2)&=\theta C_{\mathrm{cl}}+m_+h+o(h).
\end{align*}
All attainment conclusions above transfer with these parameter changes.
\end{corollary}

The rest of the paper follows the order of the main results.
Section~\ref{sec:preliminaries} establishes the estimates and
compactness arguments used throughout the proofs.
Section~\ref{sec:endpoint-proofs} proves attainment and compactness
at the finite thresholds and computes the one-sided slopes. Section~\ref{sec:strong-proofs} establishes
the two-term strong-penalty expansion, the transition as $p\downarrow4$
and their limiting profiles, followed by the anisotropic transfer.
\section{Preliminary results and compactness}\label{sec:preliminaries}

The bounded-domain Moser--Trudinger inequality~\cite{Moser} gives the following
interpolation estimate, with explicit control of the growth in $q$.

\Needspace{6\baselineskip}
\begin{lemma}\label{lem:gn}
There is an absolute constant $D\ge1$ such that, for every $q\ge r\ge2$ and $u\in H^1(\R^2)$,
\begin{equation}
 \int_{\R^2} |u|^q\dd\le(D\sqrt q)^q
 \|\nabla u\|_2^{q-r}\|u\|_r^r.
 \label{eq:gn-growth}
\end{equation}
\end{lemma}

\begin{proof}
The assertion is trivial for $u=0$. For $u\ne0$, both $\|\nabla u\|_2^2$ and
$\|u\|_r^r$ are positive. By amplitude and spatial scaling, it suffices
to consider $\|\nabla w\|_2=\|w\|_r=1$. Specifically, set
\[
 w(x)=\|\nabla u\|_2^{-1}u(bx),\qquad
 b^2=\|u\|_r^r\|\nabla u\|_2^{-r}.
\]
Let $A=\{|w|>1\}$ and $z=(|w|-1)_+$. Then $|A|\le1$,
$\|\nabla z\|_2^2\le1$, and Schwarz rearrangement puts $z$ in $H^1_0$ of the ball
of area $|A|$ without increasing its energy. The scale-invariant
bounded-domain Moser--Trudinger inequality gives
\[
 \int_A e^{4\pi z^2}\dd\le C_0|A|\le C_0.
\]
If $|A|=0$, this step is unnecessary. On $A$,
$2\pi|w|^2\le4\pi z^2+4\pi$, and therefore
$\int_Ae^{2\pi|w|^2}\dd\le C_0e^{4\pi}$.
The elementary bound
\[
 s^q\le\left(\frac q{4\pi e}\right)^{q/2}e^{2\pi s^2},
 \qquad s\ge0,
\]
and $\int_{A^c}|w|^q\dd\le \|w\|_r^r=1$ yield
\[
 \|w\|_q^q\le1+C_0e^{4\pi}
 \left(\frac q{4\pi e}\right)^{q/2}\le(D\sqrt q)^q
 \quad(q\ge r\ge2),
\]
with an absolute constant $D$. Rescaling to $u$ gives
\eqref{eq:gn-growth}.
\end{proof}

\Needspace{6\baselineskip}
\begin{lemma}\label{lem:small}
There are constants $\varepsilon_0>0$ and $C>0$ such that, whenever
$\|\nabla u\|_2^2\le\varepsilon_0$,
\begin{equation}
 0\le\int_{\R^2}(e^{4\pi u^2}-1-4\pi u^2)\dd
 \le C\|\nabla u\|_2^{4-r}\|u\|_r^r,\qquad 2\le r\le4.
 \label{eq:small-r}
\end{equation}
\begin{equation}
 0\le\int_{\R^2}(e^{4\pi u^2}-1-4\pi u^2-8\pi^2u^4)\dd
 \le C\|\nabla u\|_2^2\|u\|_4^4.
 \label{eq:small-four}
\end{equation}
\begin{equation}
 \begin{aligned}
 0&\le\int_{\R^2}\left(e^{4\pi u^2}-1-4\pi u^2-8\pi^2u^4
                  -\frac{32\pi^3}{3}u^6\right)\dd\\
  &\le C\|\nabla u\|_2^2\|u\|_6^6.
 \end{aligned}
 \label{eq:small-six}
\end{equation}
The constants may be taken independent of $r\in[2,4]$.
\end{lemma}

\begin{proof}
Put $g=\|\nabla u\|_2^2$. By monotone convergence,
\[
 \int_{\R^2}(e^{4\pi u^2}-1-4\pi u^2)\dd=\sum_{j=2}^{\infty}\frac{(4\pi)^j}{j!}\|u\|_{2j}^{2j}.
\]
\reftext{Lemma}{lem:gn} and $j!\ge(j/e)^j$ give, with
$K=8\pi eD^2$,
\[
 \frac{(4\pi)^j}{j!}\|u\|_{2j}^{2j}
 \le K^j g^{j-r/2}\|u\|_r^r,\qquad 2j\ge r.
\]
Choose $\varepsilon_0$ so that $K\varepsilon_0\le1/2$.
Summing from $j=2$ gives
$\int_{\R^2}(e^{4\pi u^2}-1-4\pi u^2)\dd\le2K^2g^{2-r/2}\|u\|_r^r$.
Taking $r=4$ and summing from $j=3$ instead gives
\eqref{eq:small-four}. Taking $r=6$ and summing from $j=4$
gives \eqref{eq:small-six}, with the same choice of sufficiently
small $\varepsilon_0$. If $g=0$, the function is zero and the
assertions follow directly.
\end{proof}

Ruf's inequality~\cite{Ruf} bounds the exponential integral
uniformly on $\A$, and Sobolev embedding bounds every finite power
moment. As a supremum of affine functions, $S(\cdot,p)$ is finite,
convex and nonincreasing. Moreover,
\begin{equation}
 |S(\lambda,p)-S(\mu,p)|
 \le\Bigl(\sup_{u\in\A}\|u\|_p^p\Bigr)|\lambda-\mu|.
 \label{eq:lipschitz}
\end{equation}

\Needspace{6\baselineskip}
\begin{lemma}\label{lem:normalization}
For $p>2$, one has $S(\lambda,p)\ge4\pi$ for every real $\lambda$.
For every $p\ge2$, any extremal with positive value satisfies $\|u\|_{H^1}=1$. Its Schwarz rearrangement is also an extremal and has the same gradient energy and nonlinear moments. Whenever $S(\lambda,p)>0$, there is a centred
radial maximizing sequence with full norm.
\end{lemma}

\begin{proof}
Fix $v\in H^1(\R^2)$ with $\|v\|_2=1$. For $t>0$, define
\begin{equation}
 U_{t,v}(x)=\frac{\sqrt t\,v(\sqrt t\,x)}{\sqrt{1+t\|\nabla v\|_2^2}}.
 \label{eq:test-family}
\end{equation}
This family has full norm, and
\begin{equation}
 \begin{aligned}
 \|\nabla U_{t,v}\|_2^2&=\frac{t\|\nabla v\|_2^2}{1+t\|\nabla v\|_2^2},\\
 \|U_{t,v}\|_2^2&=\frac1{1+t\|\nabla v\|_2^2},\\
 \|U_{t,v}\|_r^r&=\frac{t^{r/2-1}\|v\|_r^r}{(1+t\|\nabla v\|_2^2)^{r/2}}.
 \end{aligned}
 \label{eq:test-moments}
\end{equation}
For $p>2$, \reftext{Lemma}{lem:small} gives $\int_{\R^2}(e^{4\pi U_{t,v}^2}-1)\dd\to4\pi$ and
$\|U_{t,v}\|_p^p\to0$ as $t\downarrow0$. Hence $S(\lambda,p)\ge4\pi$.

Suppose that $u\ne0$ is admissible, $J_{\lambda,p}(u)>0$, and
$\|u\|_{H^1}<1$. Set
$a^2=(1-\|\nabla u\|_2^2)/\|u\|_2^2>1$ and $w(x)=u(x/a)$.
Then $\|\nabla w\|_2^2=\|\nabla u\|_2^2$, $\|w\|_2^2=a^2\|u\|_2^2=1-\|\nabla u\|_2^2$, and
\begin{equation}
 J_{\lambda,p}(w)=a^2J_{\lambda,p}(u)>J_{\lambda,p}(u).
 \label{eq:spatial-normalization}
\end{equation}
Thus every extremal of positive value has full norm.
Schwarz rearrangement preserves the exponential integral and every
$L^q$ norm, and does not increase the gradient norm. Rearrange a maximizing sequence of positive value, then dilate
spatially to restore full norm. This gives the required maximizing
sequence. If $u$ is an extremal, its rearrangement has the same positive
value and is admissible. Both functions therefore have full norm;
equality of their masses gives equality of their gradient energies.
\end{proof}

The quartic threshold estimate and the unit-mass problem both use the radial function
\begin{equation}
 v_0(x)=\sqrt{\frac3\pi}\,(1+|x|^2)^{-2}.
 \label{eq:rational-test}
\end{equation}
Direct integration gives
\begin{equation}
 \|v_0\|_2=1,\quad \|\nabla v_0\|_2^2=\frac{12}{5},\quad
 \|v_0\|_4^4=\frac9{7\pi},\quad
 8\pi^2\|v_0\|_4^4-4\pi \|\nabla v_0\|_2^2=\frac{24\pi}{35}>0.
 \label{eq:rational-moments}
\end{equation}

\begin{figure}[!htbp]
\centering
\begin{minipage}[t]{0.485\textwidth}
\centering
\begin{tikzpicture}
\begin{axis}[
 width=0.97\linewidth,height=5.4cm,
 xmin=0,xmax=10,ymin=0,ymax=0.58,
 axis lines=left,axis line style={black},
 xlabel={$r$},ylabel={$U_{\delta^2,v_0}(r)$},
 title={(a) Normalized spatial spreading},
 tick label style={font=\footnotesize},label style={font=\footnotesize},
 title style={font=\footnotesize},
 xtick={0,2,4,6,8,10},ytick={0,0.1,0.2,0.3,0.4,0.5},
 yticklabel style={/pgf/number format/fixed},
 legend style={font=\footnotesize,draw=none,fill=none,
               at={(0.97,0.97)},anchor=north east},
 samples=161,domain=0:10,clip=true]
 \addplot[black,line width=1.1pt] {sqrt(3/pi)*(1+x^2)^(-2)/sqrt(1+12/5)};
 \addlegendentry{$\delta=1$}
 \addplot[color={rgb,255:red,35;green,83;blue,125},
          line width=1.2pt,dashed]
          {0.5*sqrt(3/pi)*(1+0.5^2*x^2)^(-2)/sqrt(1+(12/5)*0.5^2)};
 \addlegendentry{$\delta=\tfrac12$}
 \addplot[color={rgb,255:red,163;green,92;blue,35},
          line width=1.2pt,dashdotted]
          {0.25*sqrt(3/pi)*(1+0.25^2*x^2)^(-2)/sqrt(1+(12/5)*0.25^2)};
 \addlegendentry{$\delta=\tfrac14$}
\end{axis}
\end{tikzpicture}
\end{minipage}\hfill
\begin{minipage}[t]{0.485\textwidth}
\centering
\begin{tikzpicture}
\begin{axis}[
 width=0.97\linewidth,height=5.4cm,
 xmin=0,xmax=4,ymin=0,ymax=1.05,
 axis lines=left,axis line style={black},
 xlabel={$r$},ylabel={rescaled amplitude},
 title={(b) Exact rescaling},
 tick label style={font=\footnotesize},label style={font=\footnotesize},
 title style={font=\footnotesize},
 xtick={0,1,2,3,4},ytick={0,0.5,1},
 yticklabel style={/pgf/number format/fixed},
 legend style={font=\footnotesize,draw=none,fill=none,
               at={(0.97,0.97)},anchor=north east},
 samples=161,domain=0:4,clip=true]
 \addplot[black,line width=1.2pt] {sqrt(3/pi)*(1+x^2)^(-2)};
 \addlegendentry{$v_0(r),\quad\delta=1,\tfrac12,\tfrac14$}
\end{axis}
\end{tikzpicture}
\end{minipage}
\caption{Explicit rational tests. Decreasing $\delta$ lowers and
broadens the profile; exact rescaling maps all three curves to $v_0$.}
\label{fig:normalized-spreading}
\end{figure}
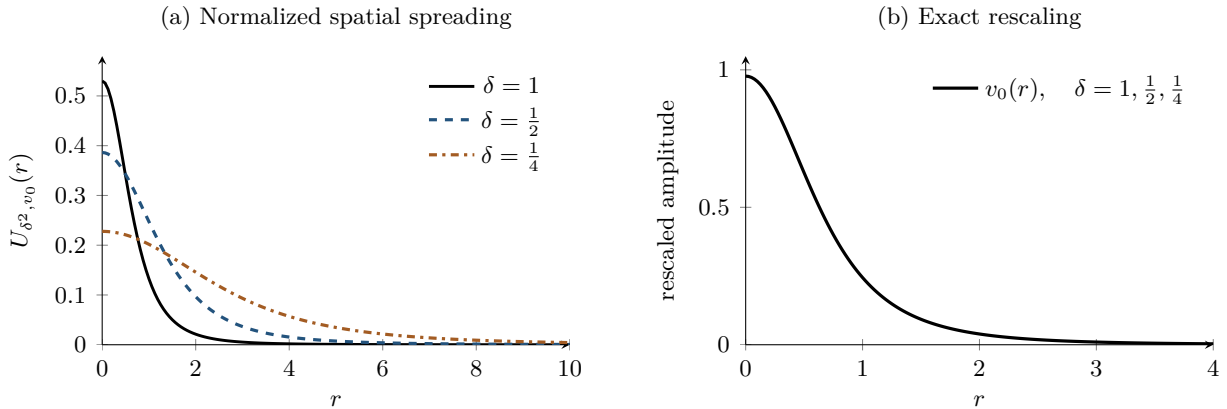

For the rational test $v_0$, writing $r=|x|$ and $t=\delta^2$ gives
\[
 U_{\delta^2,v_0}(r)=
 \frac{\delta\sqrt{3/\pi}}{\sqrt{1+(12/5)\delta^2}}
 (1+\delta^2r^2)^{-2},\qquad
 \sqrt{1+(12/5)\delta^2}\,\delta^{-1}
 U_{\delta^2,v_0}(r/\delta)=v_0(r).
\]
\FloatBarrier
Figure~\ref{fig:normalized-spreading} illustrates this exact rescaling.

\FloatBarrier
Radial monotonicity now gives the tail control needed for compactness.

For a centred radial function with $\|u\|_2^2\le1$,
\begin{equation}
 u(r)^2\le\frac1{\pi r^2},\qquad r>0.
 \label{eq:radial-bound}
\end{equation}
Consequently, for a bounded sequence of such functions converging weakly
in $H^1$ to $u$, local Rellich compactness and the tail estimate
\[
 \int_{|x|>R}|u_n|^q\dd
 \le(\pi R^2)^{-(q-2)/2}\|u_n\|_2^2,\qquad q>2,
\]
give strong convergence in every finite $L^q$, $q>2$, after passage to
a subsequence.
We also use the following integrability property for a fixed function:
\begin{equation}
 \int_{\R^2}\bigl(e^{a w^2}-1\bigr)\dd<\infty
 \quad\text{for each fixed }w\in H^1(\R^2)\text{ and }a>0.
 \label{eq:fixed-integrability}
\end{equation}
To see this, choose a truncation $z=(|w|-L)_+$ with sufficiently small
$H^1$ norm. On $\{|w|>L\}$, use
$w^2\le2z^2+2L^2$ and Ruf's inequality for the normalized $z$;
the set has finite measure. On its complement,
$e^{aw^2}-1\le C_{a,L}w^2$.

\Needspace{6\baselineskip}
\begin{lemma}\label{lem:remainder-continuity}
Let $u_n$ be centred radial functions with $\|u_n\|_{H^1}=1$, and
let $u_n\weak u$ in $H^1$.
\begin{enumerate}[label=\textup{(\roman*)}]
\item If $u\ne0$, then
\[
 \int_{\R^2}(e^{4\pi u_n^2}-1-4\pi u_n^2)\dd
 \longrightarrow\int_{\R^2}(e^{4\pi u^2}-1-4\pi u^2)\dd.
\]
\item If $u=0$ and $\limsup_n\|\nabla u_n\|_2^2<1$, then
\[
 \int_{\R^2}(e^{4\pi u_n^2}-1-4\pi u_n^2)\dd\longrightarrow0.
\]
\end{enumerate}
\end{lemma}

\begin{proof}
Outside a fixed sufficiently large ball, \eqref{eq:radial-bound} and
$0\le e^{4\pi s^2}-1-4\pi s^2\le Cs^4$ for $|s|\le1$ yield
\begin{equation}
 \int_{|x|>R}\bigl(e^{4\pi u_n^2}-1-4\pi u_n^2\bigr)\dd
 \le\frac C{R^2}.
 \label{eq:remainder-tail}
\end{equation}
The same bound applies to $u$.

For \hyperref[lem:remainder-continuity]{\textup{(i)}}, set $s=\|u\|_{H^1}^2>0$ and $z_n=u_n-u$. The Hilbert-space
identity gives $\|z_n\|_{H^1}^2\to1-s<1$. Choose $q>1$,
$\eta>0$, and a H\"older exponent $a>1$ such that
$qa(1+\eta)(1-s)<1$. On each fixed ball,
\[
 e^{4\pi q u_n^2}
 \le e^{4\pi q(1+\eta)z_n^2}
 e^{4\pi q(1+\eta^{-1})u^2}.
\]
H\"older's inequality, Ruf's estimate for $z_n$, and
\eqref{eq:fixed-integrability} give a uniform integral bound for the
left side. Thus the nonlinear remainders are uniformly integrable on
the ball and converge there in $L^1$. Combine this with
\eqref{eq:remainder-tail}.

For \hyperref[lem:remainder-continuity]{\textup{(ii)}}, fix $R>0$. Local $L^2$ convergence and radial
monotonicity imply $u_n(R)\to0$. The function
$z_n=(u_n-u_n(R))_+$ on $B_R$ has zero trace and
$\int_{B_R}|\nabla z_n|^2\dd\le \|\nabla u_n\|_2^2\le a_0<1$ for large $n$. Choose $q>1$
and $\eta>0$ so that $q(1+\eta)a_0<1$, and apply the
bounded-domain Moser--Trudinger inequality to
\[
 e^{4\pi q u_n^2}
 \le e^{4\pi q(1+\eta)z_n^2}
 e^{4\pi q(1+\eta^{-1})u_n(R)^2}
 \quad\text{on }B_R.
\]
Uniform integrability and local convergence to zero give convergence
of the remainder integral on $B_R$ to zero. Then let $R\to\infty$
in \eqref{eq:remainder-tail}.
\end{proof}

\Needspace{6\baselineskip}
\begin{lemma}\label{lem:weakzero}
Suppose that $u_n$ are centred radial, $\|u_n\|_{H^1}=1$,
$u_n\weak0$ in $H^1$, and
\begin{equation}
 \liminf_n \int_{\R^2}(e^{4\pi u_n^2}-1)\dd\ge4\pi.
 \label{eq:weakzero-level}
\end{equation}
Then $\|\nabla u_n\|_2^2\to0$ and $\|u_n\|_2^2\to1$.
\end{lemma}

\begin{proof}
Pass to a subsequence with $\|\nabla u_n\|_2^2\to g\in[0,1]$.
If $g<1$, Lemma~\ref{lem:remainder-continuity}\textup{(ii)} gives
$\int_{\R^2}(e^{4\pi u_n^2}-1)\dd\to4\pi(1-g)$. Condition \eqref{eq:weakzero-level} forces
$g=0$.

Suppose $g=1$, so $\|u_n\|_2^2\to0$. If the sequence concentrates at the origin, in the sense that
\[
 \int_{|x|>R}\bigl(|\nabla u_n|^2+u_n^2\bigr)\,dx\longrightarrow0
 \qquad\text{for every }R>0,
\]
Nguyen's estimate \eqref{eq:concentration} gives
$\limsup \int_{\R^2}(e^{4\pi u_n^2}-1)\dd\le C_{\mathrm{cl}}<4\pi$, a contradiction.
If it does not concentrate, there are a fixed $R>0$, an $\varepsilon>0$,
and a subsequence with
\[
 \int_{|x|>R}|\nabla u_n|^2\dd\ge\varepsilon.
\]
The truncated functions on $B_R$ then have gradient energy at most
$1-\varepsilon$, and $u_n(R)\to0$. The local uniform-integrability
argument in \reftext{Lemma}{lem:remainder-continuity} shows that
$\int_{B_R}(e^{4\pi u_n^2}-1)\dd\to0$.
Outside $B_R$, radial monotonicity gives
$\sup u_n^2\le \|u_n\|_2^2/(\pi R^2)\to0$; hence
\[
 \int_{|x|>R}(e^{4\pi u_n^2}-1)\dd\le C \|u_n\|_2^2\to0.
\]
This contradicts \eqref{eq:weakzero-level}. Every subsequential limit
of $\|\nabla u_n\|_2^2$ is therefore zero, proving the assertion.
\end{proof}

\Needspace{6\baselineskip}
\begin{lemma}\label{lem:nonzero-compactness}
Let $p>2$ and $\lambda_n\to\lambda\in\R$, and let $u_n$ be centred radial
functions with full norm such that
\[
 J_{\lambda_n,p}(u_n)-S(\lambda_n,p)\longrightarrow0.
\]
If $u_n\weak u\ne0$ in $H^1$, then $u_n\to u$ strongly in $H^1$
and $u$ is an extremal at $\lambda$.
\end{lemma}

\begin{proof}
Take a subsequence with $\|u_n\|_2^2\to m$. \reftext{Lemma}{lem:remainder-continuity},
strong $L^p$ convergence, and \eqref{eq:lipschitz} imply
\begin{equation}
 S(\lambda,p)=4\pi m+
 \int_{\R^2}(e^{4\pi u^2}-1-4\pi u^2-\lambda|u|^p)\dd.
 \label{eq:mass-splitting}
\end{equation}
Weak lower semicontinuity gives $\|u\|_2^2\le m$ and $\|\nabla u\|_2^2\le1-m$.
Since $u\ne0$, $\|\nabla u\|_2^2>0$, so $m<1$. Consequently
\begin{equation}
 \int_{\R^2}(e^{4\pi u^2}-1-4\pi u^2-\lambda|u|^p)\dd
 =S(\lambda,p)-4\pi m\ge4\pi(1-m)>0.
 \label{eq:positive-nonlinear}
\end{equation}
Set $\tau=m/\|u\|_2^2\ge1$ and $w(x)=u(x/\sqrt\tau)$.
Then $w\in\A$, $\|w\|_2^2=m$, and
\[
 J_{\lambda,p}(w)=4\pi m+
 \tau\int_{\R^2}(e^{4\pi u^2}-1-4\pi u^2-\lambda|u|^p)\dd.
\]
By \eqref{eq:mass-splitting} and \eqref{eq:positive-nonlinear},
$J_{\lambda,p}(w)>S(\lambda,p)$ if $\tau>1$.
Thus $\tau=1$ and $\|u\|_2^2=m$.
It follows that $J_{\lambda,p}(u)=S(\lambda,p)$. Its positive value
and \reftext{Lemma}{lem:normalization} force $\|u\|_{H^1}=1$.
The $L^2$ and gradient norms therefore converge to those of $u$,
which proves strong $H^1$ convergence.
\end{proof}

For $p=2$, the tests \eqref{eq:test-family} and the sharp concentrating
sequences for \eqref{eq:concentration} give
\[
 S(\lambda,2)\ge\max\{4\pi-\lambda,C_{\mathrm{cl}}\}.
\]
The mass of a concentrating sequence tends to zero, so its fixed quadratic
penalty does not change the limiting value. The following compactness
statement applies above these two barriers.

\Needspace{6\baselineskip}
\begin{lemma}\label{lem:quadratic-compactness}
For the quadratic problem, a normalized centred radial maximizing sequence
with a nonzero weak limit converges strongly in $H^1$ along a subsequence.
In particular, if
\begin{equation}
 S(\lambda,2)>\max\{C_{\mathrm{cl}},4\pi-\lambda\},
 \label{eq:quadratic-barriers}
\end{equation}
then the supremum is attained, and every such maximizing sequence is
precompact in $H^1$.
\end{lemma}

\begin{proof}
Set $d=4\pi-\lambda$. For a nonzero weak limit $u$, pass to a subsequence
with $\|u_n\|_2^2\to m$. \reftext{Lemma}{lem:remainder-continuity} gives
$S(\lambda,2)=dm+\int_{\R^2}(e^{4\pi u^2}-1-4\pi u^2)\dd$, while $\|u\|_2^2\le m$ and $\|\nabla u\|_2^2\le1-m$.
Let $\tau=m/\|u\|_2^2$ and $w(x)=u(x/\sqrt\tau)$. Then $w\in\A$ and
\[
 J_{\lambda,2}(w)=dm+\tau\int_{\R^2}(e^{4\pi u^2}-1-4\pi u^2)\dd.
\]
Since $u\ne0$, the remainder integral is positive, so comparison with
the supremum gives $\tau=1$.
The limit is an extremal. Its positive value forces full norm by
\reftext{Lemma}{lem:normalization}, and both Hilbert-space norms converge.
The same argument applies when the parameters converge to a finite limit, since $S(\cdot,2)$ is $1$-Lipschitz and the $L^2$ masses are bounded by one.

If the weak limit is zero, pass to a subsequence with $\|\nabla u_n\|_2^2\to g$.
When $g<1$, \reftext{Lemma}{lem:remainder-continuity} gives
$J_{\lambda,2}(u_n)\to d(1-g)\le\max\{d,C_{\mathrm{cl}}\}$.
When $g=1$, the $L^2$ mass tends to zero. The concentrating alternative has
$\limsup \int_{\R^2}(e^{4\pi u_n^2}-1)\dd\le C_{\mathrm{cl}}$; the nonconcentrating alternative has
$\int_{\R^2}(e^{4\pi u_n^2}-1)\dd\to0$, by the truncation argument in \reftext{Lemma}{lem:weakzero}.
Thus a weak zero limit cannot have a value exceeding both barriers.
\end{proof}

The endpoint proofs also require the following strict estimates.

Recall $B_2$ and $\lambda_{\mathrm{loc}}$ from \eqref{eq:gn-sharp}.
The estimate in \reftext{Lemma}{lem:gn} shows that $B_2$ is finite.
The rational test \eqref{eq:rational-test} gives
\begin{equation}
 B_2\ge\frac{15}{28\pi}>\frac1{2\pi}.
 \label{eq:gn-rational-bound}
\end{equation}

\Needspace{6\baselineskip}
\begin{lemma}\label{lem:gn-moment}
There is a centred radial optimizer $V$ for $B_2$ with
$\|\nabla V\|_2=\|V\|_2=1$. It satisfies
\begin{equation}
 \|V\|_4^4=B_2,\qquad \|V\|_6^6\ge\frac43 B_2^2.
 \label{eq:gn-sixth}
\end{equation}
\end{lemma}

\begin{proof}
Normalize a maximizing sequence by amplitude and spatial scaling to have
unit $L^2$ mass and unit gradient norm. Rearrange it and use a
mass-preserving dilation to restore the unit gradient norm.
These operations do not decrease the quotient in \eqref{eq:gn-sharp},
which remains bounded above by $B_2$. The resulting sequence is therefore
still maximizing, with unit $L^2$ mass and unit gradient norm. Radial compactness gives a weak
limit $V$ and strong $L^4$ convergence, so
\[
 B_2=\|V\|_4^4\le B_2 \|\nabla V\|_2^2\|V\|_2^2\le B_2.
\]
Both weak lower semicontinuity inequalities are equalities. Hence the
sequence converges strongly in $H^1$, and $\|\nabla V\|_2=\|V\|_2=1$.

Differentiating the quotient at $V$ yields
$-\Delta V+V=(2/B_2)V^3$. Standard local elliptic regularity and positivity
give a smooth positive radial solution
$Q=\sqrt{2/B_2}\,V$ of
\begin{equation}
 Q''+r^{-1}Q'=Q-Q^3.
 \label{eq:gn-ode}
\end{equation}
The function $Q$ is nonincreasing and tends to zero. In fact $Q'<0$ for
$r>0$: if $Q'$ vanished at an interior point, its nonpositivity would
also give $Q''=0$ there, and \eqref{eq:gn-ode} would give $Q=1$.
Uniqueness for the ordinary differential equation would force $Q\equiv1$.
Since $\|Q\|_4^4=2\|Q\|_2^2$, one also has $Q(0)^2>2$.

Set $z=-Q'/(rQ)>0$. Equation~\eqref{eq:gn-ode} becomes
\begin{equation}
 rz'=r^2z^2-2z+Q^2-1,
 \qquad (Q^2)'=-2rzQ^2.
 \label{eq:riccati}
\end{equation}
Writing $a=Q(0)$, the regular expansion at the origin gives
\[
 z(r)=\frac{a^2-1}{2}-\frac{a^4-1}{16}r^2+O(r^4).
\]
Thus $z'<0$ initially. We claim it remains negative. At any point where
$z'=0$, differentiation of \eqref{eq:riccati} gives
$z''=2z(z-Q^2)$. At a hypothetical first zero $r_0$ of $z'$, one has $z''(r_0)\ge0$,
hence $z(r_0)\ge Q(r_0)^2$.
If $z=Q^2$ there, a second differentiation gives
$z'''(r_0)=4r_0Q(r_0)^6>0$, which contradicts $z'<0$ immediately to
the left. If $z>Q^2$, then $z'>0$ immediately to the right. As long as
$z'>0$, the difference $z-Q^2$ increases; a subsequent zero of $z'$ would
have $z''>0$ and could not be a crossing from positive values. Hence
$z'>0$ thereafter. With $z$ bounded below by a positive constant,
\eqref{eq:riccati} gives $z'\ge rz^2/2$ for all sufficiently large $r$.
Integrating the inequality for $1/z$ forces a finite-radius blow-up,
contrary to smoothness and positivity of $Q$. This proves $z'<0$.

Let $s=r^2$, $f(s)=Q(\sqrt{s})^2$, and $b=f(0)$.
The change of variable $t=\log(b/f(s))$ has
$dt/ds=z(\sqrt{s})>0$. Its inverse derivative
$w(t)=ds/dt$ is increasing. For $k=1,2,3$ put
$I_k=\int_0^\infty f(s)^k\,ds$ and define the positive measure
$d\nu=w(0)\delta_0+dw$ on $[0,\infty)$. Integration by parts yields
\begin{equation}
 \frac{kI_k}{b^k}=\int_0^\infty e^{-kt}\,d\nu(t).
 \label{eq:moment-laplace}
\end{equation}
The boundary term at infinity vanishes: finiteness of $I_k$ and
monotonicity of $w$ imply $e^{-kt}w(t)\to0$.
Cauchy--Schwarz applied to \eqref{eq:moment-laplace} gives
$3I_1I_3\ge4I_2^2$. Since planar radial integration is $\pi\,ds$,
\[
 \|Q\|_2^2\|Q\|_6^6\ge\frac43\|Q\|_4^8.
\]
Amplitude normalization gives \eqref{eq:gn-sixth}.
\end{proof}

For $2<p\le4$, let $\lamstar\in(0,\infty)$ be the threshold of
Chen--Jiang--Lu--Zhu~\cite[Theorem 1.1]{Chen}. Above this threshold
they proved that $S(\lambda,p)=4\pi$ and that the supremum is not
attained. Continuity gives $S(\lamstar,p)=4\pi$. To identify the value
below the threshold, take $\lambda<\mu<\lamstar$ and an extremal $u_\mu$
at $\mu$. Lemma~\ref{lem:normalization} gives $S(\mu,p)\ge4\pi$, so
\begin{equation}
 S(\lambda,p)\ge S(\mu,p)+(\mu-\lambda)\|u_\mu\|_p^p>4\pi.
 \label{eq:attained-comparison}
\end{equation}
Consequently,
\begin{equation}
 \lambda_*(p)=\sup\{\lambda\in\R:S(\lambda,p)>4\pi\}.
 \label{eq:threshold-level-definition}
\end{equation}
For $p>4$, the same theorem gives an extremal at every finite parameter.
Taking $\mu=\lambda+1$ in \eqref{eq:attained-comparison} gives
$S(\lambda,p)>4\pi$ for every $\lambda\in\R$.

\begin{proposition}\label{prop:quartic-gap}
The quartic threshold satisfies
\begin{equation}
 \lambda_*(4)>\lambda_{\mathrm{loc}}>0.
 \label{eq:quartic-strict}
\end{equation}
\end{proposition}

\begin{proof}
Use $V$ from \reftext{Lemma}{lem:gn-moment} in the normalized test family
\eqref{eq:test-family}. The expansion through degree six gives
\begin{equation}
 J_{\lambda_{\mathrm{loc}},4}(U_{t,V})
 =4\pi+\left(\frac{32\pi^3}{3}\|V\|_6^6-4\pi\right)t^2+O(t^3).
 \label{eq:quartic-test-second}
\end{equation}
For this fixed profile the remainder follows from \reftext{Lemma}{lem:gn}
by summing the exponential series from degree eight.
Equations~\eqref{eq:gn-rational-bound} and \eqref{eq:gn-sixth} give
\[
 \frac{32\pi^3}{3}\|V\|_6^6-4\pi
 \ge4\pi\left(\frac{32\pi^2}{9}B_2^2-1\right)
 \ge\frac{4\pi}{49}>0.
\]
Hence $S(\lambda_{\mathrm{loc}},4)>4\pi$. By continuity, the value
remains above $4\pi$ at some parameter larger than $\lambda_{\mathrm{loc}}$.
The definition of $\lambda_*(4)$ gives \eqref{eq:quartic-strict}.
\end{proof}

The threshold quotients may be restricted to functions of full norm:
spatial dilation to full norm increases both quotients in
\eqref{eq:quadratic-quotients}.

\Needspace{6\baselineskip}
\begin{lemma}
\label{lem:quadratic-explicit-gap}
The quadratic threshold quotient satisfies
\begin{equation}
 D_*>\frac{4743\pi}{1000}>\frac{128\pi}{27}\ge8\pi^2B_2.
 \label{eq:quadratic-explicit-gap}
\end{equation}
\end{lemma}

\begin{proof}
We first record an explicit upper bound for $B_2$. Use the unitary Fourier
transform $\mathscr F u(\xi)=(2\pi)^{-1}\int_{\R^2} e^{-ix\cdot\xi}u(x)\,dx$.
The sharp Hausdorff--Young inequality of Beckner~\cite{Beckner}, followed
by H\"older's inequality, gives, for every $a>0$,
\[
 \|u\|_4^4\le\frac{4}{27\pi^2}\|\mathscr F u\|_{4/3}^4
 \le\frac{4}{27\pi a}\bigl(\|\nabla u\|_2^2+a\|u\|_2^2\bigr)^2.
\]
Here $\int_{\R^2}(a+|\xi|^2)^{-2}\,d\xi=\pi/a$.
Taking $a=\|\nabla u\|_2^2/\|u\|_2^2$ yields
\begin{equation}
 B_2\le\frac{16}{27\pi}.
 \label{eq:gn-explicit-upper}
\end{equation}

Define $\varphi(r)=f((1+r)^{-1})$, where
\begin{align}
 100f(\xi)={}&-480\xi^3+10710\xi^4-49644\xi^5
 +109620\xi^6-130200\xi^7\nonumber\\
 &\hspace{8mm}+80235\xi^8-20030\xi^9-111\xi^{10}.
 \label{eq:quadratic-rational-profile}
\end{align}
The radial function $\varphi$ belongs to $H^1(\R^2)$.
The change of variables $\xi=(1+r)^{-1}$ gives
\begin{align}
 e_0:=\frac{\|\nabla \varphi\|_2^2}{\pi}
 &=2\int_0^1\xi(1-\xi)f'(\xi)^2\,d\xi
 =\frac{4087991}{9509500},\nonumber\\
 m_0:=\frac{\|\varphi\|_2^2}{\pi}
 &=2\int_0^1\frac{1-\xi}{\xi^3}f(\xi)^2\,d\xi
 =\frac{787647607}{2282280000},\nonumber\\
 h_j:=\frac{\|\varphi\|_{2j}^{2j}}{\pi}
 &=2\int_0^1\frac{1-\xi}{\xi^3}f(\xi)^{2j}\,d\xi.
 \label{eq:quadratic-test-integrals}
\end{align}
Set
\[
 u(x)=\sqrt{\frac{3}{4\pi e_0}}\,\varphi(b|x|),
 \qquad b^2=\frac{3m_0}{e_0}.
\]
Then $\|\nabla u\|_2^2=3/4$ and $\|u\|_2^2=1/4$. Truncating the exponential series
after degree thirty-six gives
\begin{equation}
 \frac{\int_{\R^2}(e^{4\pi u^2}-1-4\pi u^2)\dd}{\pi \|\nabla u\|_2^2}
 \ge \frac4{m_0}\sum_{j=2}^{18}
 \frac{3^{j-2}h_j}{j!e_0^{j-1}}
 >\frac{4743}{1000}.
 \label{eq:quadratic-finite-certificate}
\end{equation}
Set $a_j=4\cdot3^{j-2}h_j/(m_0j!e_0^{j-1})$.
The integer bounds below certify \eqref{eq:quadratic-finite-certificate}.
\begin{center}
\begin{tabular}{c|r}
 Range of $j$ & Lower bound for $10^6\sum a_j$\\ \hline
 $2$--$4$ & $2243960$\\
 $5$--$8$ & $2023845$\\
 $9$--$12$ & $444324$\\
 $13$--$16$ & $30329$\\
 $17$--$18$ & $704$
\end{tabular}
\end{center}
Their sum is $4743162>4743000$. Every entry uses rational arithmetic:
if $f(\xi)^{2j}=\sum_{k=6j}^{20j}c_{j,k}\xi^k$, then
\[
 h_j=2\sum_{k=6j}^{20j}\frac{c_{j,k}}{(k-2)(k-1)}.
\]
This formula reduces every table entry to a finite rational calculation
and verifies \eqref{eq:quadratic-finite-certificate}. Combining it with
\eqref{eq:gn-explicit-upper} proves \eqref{eq:quadratic-explicit-gap}.
\end{proof}

Let $G$ be the fundamental solution of $-\Delta+1$ on $\R^2$.
In terms of the modified Bessel function $K_0$ of the second kind,
\begin{equation}
 G(r)=\frac1{2\pi}K_0(r),\qquad
 G(r)=-\frac1{2\pi}\log r+A_0+O(r^2|\log r|),\qquad
 A_0=\frac{\log2-\gamma_{\mathrm E}}{2\pi}.
 \label{eq:green-normalization}
\end{equation}
Thus $(-\Delta+1)G=\delta_0$, which fixes both the flux and the regular part. This is the normalization used in \eqref{eq:concentration}. With the unitary Fourier transform used above,
$\mathscr F G(\xi)=(2\pi)^{-1}(1+|\xi|^2)^{-1}$. Parseval's identity gives
\begin{equation}
 \int_{\R^2} G^2\,dx
 =\int_{\R^2}|\mathscr F G(\xi)|^2\,d\xi
 =(2\pi)^{-2}\int_{\R^2}(1+|\xi|^2)^{-2}\,d\xi=\frac1{4\pi}.
 \label{eq:green-mass}
\end{equation}
The fourth moment of $G$ is finite and positive.

\begin{proposition}
\label{prop:quadratic-green-test}
One has
\begin{equation}
 S(4\pi+C_{\mathrm{cl}},2)>C_{\mathrm{cl}},
 \qquad \lambda_+>4\pi+C_{\mathrm{cl}}.
 \label{eq:quadratic-strict-test}
\end{equation}
If a sequence of extremals were to concentrate as in
Lemma~\ref{lem:quadratic-concentration-parameter}, its parameters would
necessarily converge to $4\pi+C_{\mathrm{cl}}$.
\end{proposition}

\begin{proof}
Let $\varepsilon\downarrow0$, $L=(\log(1/\varepsilon))^3$,
$\rho=\varepsilon L$, and $h=\log(1+\pi L^2)$. Define
\begin{equation}
 W_\varepsilon(r)=
 \begin{cases}
 G(\rho)+[h-\log(1+\pi r^2/\varepsilon^2)]/(4\pi),&r\le\rho,\\
 G(r),&r>\rho,
 \end{cases}
 \quad c^2=\|W_\varepsilon\|_{H^1}^2,\quad
 u_\varepsilon=W_\varepsilon/c.
 \label{eq:green-test}
\end{equation}
The functions $u_\varepsilon$ belong to $H^1(\R^2)$ and have full norm.
Radial integration and the exterior Green energy identity give
\begin{align}
 \int_{B_\rho}|\nabla W_\varepsilon|^2\dd
 &=\frac{h-1+(1+\pi L^2)^{-1}}{4\pi},\nonumber\\
 \int_{B_\rho^c}(|\nabla G|^2+G^2)\dd
 &=-2\pi\rho G(\rho)G'(\rho)
 =G(\rho)+O(\rho^2\log^2\rho),\nonumber\\
 \int_{B_\rho}W_\varepsilon^2\dd&=O(\rho^2\log^2\varepsilon).
 \label{eq:green-norm}
\end{align}
Put $a=4\pi c^2$ and
$\ell=a-h-4\pi G(\rho)$. Equations~\eqref{eq:green-norm} give
\begin{equation}
 a=-2\log\varepsilon+4\pi A_0+\log\pi-1+O(L^{-2})
 +O(\rho^2\log^2\varepsilon),\qquad
 \ell=-1+O(L^{-2})+O(\rho^2\log^2\varepsilon).
 \label{eq:green-matching}
\end{equation}
The errors in \eqref{eq:green-matching} are $o(a^{-2})$. Moreover,
\begin{equation}
 \|u_\varepsilon\|_2^2=\frac1a+o(a^{-2}).
 \label{eq:green-test-mass}
\end{equation}

In the inner ball set $y=\log(1+\pi r^2/\varepsilon^2)$.
The radial area element becomes $2\pi r\,dr=\varepsilon^2e^y\,dy$, and
$4\pi u_\varepsilon^2=a-2(y+\ell)+(y+\ell)^2/a$.
Since $0\le y\le h=O(\log a)$, Taylor expansion in this last square,
followed by integration against $e^{-y}$, yields
\begin{equation}
 \int_{B_\rho}(e^{4\pi u_\varepsilon^2}-1)\,dx
 =C_{\mathrm{cl}}\left(1+\frac1a+\frac9{2a^2}\right)+o(a^{-2}).
 \label{eq:green-core-expansion}
\end{equation}
Here $\varepsilon^2e^{a-2\ell}=C_{\mathrm{cl}}+o(a^{-2})$,
$\int_0^\infty e^{-y}(y-1)^2\,dy=1$, and
$\int_0^\infty e^{-y}(y-1)^4\,dy=9$.
The discarded tails are $O(L^{-2}(\log L)^4)=o(a^{-2})$;
the exponential Taylor remainder has integral $O(a^{-3})$.
The subtracted area $\pi\rho^2$ is negligible.

In the exterior region, the first two terms of the exponential
series yield the lower bound
\[
 \int_{B_\rho^c}(e^{4\pi G^2/c^2}-1)\,dx
 \ge\frac{4\pi}{c^2}\int_{B_\rho^c}G^2\,dx
 +\frac{8\pi^2}{c^4}\int_{B_\rho^c}G^4\,dx.
\]
Using \eqref{eq:green-mass}--\eqref{eq:green-core-expansion}, for each
fixed $\lambda$ one obtains
\begin{equation}
 J_{\lambda,2}(u_\varepsilon)
 \ge C_{\mathrm{cl}}+\frac{C_{\mathrm{cl}}+4\pi-\lambda}{a}
 +\frac{\frac92C_{\mathrm{cl}}+128\pi^4\int_{\R^2} G^4\,dx}{a^2}
 +o(a^{-2}).
 \label{eq:quadratic-green-lower}
\end{equation}
At $\lambda=4\pi+C_{\mathrm{cl}}$, the first correction vanishes and the
coefficient of $a^{-2}$ is positive. Choosing $\varepsilon$ sufficiently
small proves \eqref{eq:quadratic-strict-test}. This test remains strict
on an interval to the right of $4\pi+C_{\mathrm{cl}}$. Monotonicity then gives
$S(\lambda,2)>C_{\mathrm{cl}}$ for $0\le\lambda\le4\pi+C_{\mathrm{cl}}$.
\end{proof}

We next examine concentrating quadratic extremals at the upper threshold.

\begin{proposition}
\label{prop:quadratic-energy-criterion}
For every quadratic extremal $u$ of positive value, set
\begin{equation}
 \mu=4\pi\int_{\R^2} u^2e^{4\pi u^2}\,dx-\lambda \|u\|_2^2,
 \qquad
 \mathcal E=\frac{4\pi}{\mu}\int_{\R^2} u^2e^{4\pi u^2}\,dx.
 \label{eq:quadratic-multiplier}
\end{equation}
Then $\mu>0$, and
\begin{equation}
 \mu \|u\|_2^2=S(\lambda,2),\qquad
 \mu^2(\mathcal E-1)=\lambda S(\lambda,2).
 \label{eq:quadratic-energy-identity}
\end{equation}
An extremal also satisfies
\begin{equation}
 -\Delta u+u=\frac{4\pi}{\mu}u e^{4\pi u^2}
                  -\frac{\lambda}{\mu}u.
 \label{eq:quadratic-euler}
\end{equation}
\end{proposition}

\begin{proof}
The full norm gives $\|u\|_{H^1}=1$. For $s$ near one, the functions
\[
 u_s(x)=s u(b_sx),\qquad
 b_s^2=\frac{s^2\|u\|_2^2}{1-s^2\|\nabla u\|_2^2},
\]
have full norm and $u_1=u$. Stationarity along this curve gives
\[
 0=\left.\frac{d}{ds}J_{\lambda,2}(u_s)\right|_{s=1}
 =-\frac{2\int_{\R^2}(e^{4\pi u^2}-1)\dd}{\|u\|_2^2}
 +8\pi\int_{\R^2} u^2e^{4\pi u^2}\,dx+2\lambda \|\nabla u\|_2^2.
\]
Differentiation is justified by \eqref{eq:fixed-integrability}, applied
at a fixed exponent larger than $4\pi s^2$. Rearranging gives
$\mu \|u\|_2^2=\int_{\R^2}(e^{4\pi u^2}-1)\dd-\lambda \|u\|_2^2=S(\lambda,2)>0$.
The second identity in \eqref{eq:quadratic-energy-identity} follows
directly from \eqref{eq:quadratic-multiplier}. The Lagrange multiplier
equation is \eqref{eq:quadratic-euler}, with multiplier $\mu$
identified by testing against $u$. Local differentiability follows
from \eqref{eq:fixed-integrability} and the decomposition in
Lemma~\ref{lem:remainder-continuity}. Moreover,
$ue^{4\pi u^2}\in L^2(\R^2)$: on $\{|u|\le1\}$ this follows from
$u\in L^2$, and on its complement a larger fixed exponential absorbs
the polynomial factor. Resolvent regularity gives
$u\in H^2(\R^2)\subset L^\infty(\R^2)$. Elliptic bootstrapping then
gives smoothness and justifies the truncation tests below.
\end{proof}

We use the radial correctors of Mancini--Martinazzi~\cite{ManciniMartinazzi},
recalling their equations and integral identities below. The estimates
on growing intervals and the exterior bounds in
Lemma~\ref{lem:quadratic-concentration-parameter} control the additional
whole-plane terms.

\Needspace{6\baselineskip}
\begin{lemma}
\label{lem:quadratic-correctors}
Set $t(y)=\log(1+|y|^2)$ and $\psi=-t$. Define the linearized operator
\[
 \mathcal{L}H=-\Delta H-8e^{2\psi}H.
\]
Let $S,Z$ be the regular radial
solutions, with zero value and radial derivative at the origin, of
\begin{align}
 \mathcal{L}S&=4e^{2\psi}(\psi^2+\psi),\label{eq:corrector-s}\\
 \mathcal{L}Z&=4e^{2\psi}
 \bigl(S+2S^2+4\psi S+2\psi^2S+\psi^3+\tfrac12\psi^4\bigr).
 \label{eq:corrector-z}
\end{align}
They satisfy $|S|+|Z|\le C(1+t)$. Define
\begin{align}
 A&=2S+\psi^2+2\psi,\nonumber\\
 B&=2Z+2S^2+2\psi^2S+\tfrac12\psi^4
       +6\psi S+2S+\psi^2+2\psi^3.
 \label{eq:corrector-energy-coefficients}
\end{align}
Then
\begin{equation}
 4\int_{\R^2} e^{2\psi}A\,dy=0,\qquad
 4\int_{\R^2} e^{2\psi}B\,dy=4\pi.
 \label{eq:corrector-integrals}
\end{equation}
If $\mathcal{L}H=f$ on a radial interval, $H(0)=H'(0)=0$, and
$|f|\le K(1+t^m)e^{-\alpha t}$ for $m\ge0$ and $\alpha>1$, then
\begin{equation}
 |H(r)|\le C_{\alpha,m}K(1+t(r)),
 \label{eq:corrector-linear-bound}
\end{equation}
with a constant independent of the interval's endpoint.
\end{lemma}

\begin{proof}
Two fundamental solutions of the radial homogeneous equation are
\[
 z(r)=\frac{1-r^2}{1+r^2},\qquad
 q(r)=z(r)\log r+\frac{2}{1+r^2},\qquad
 zq'-z'q=\frac1r.
\]
Variation of constants gives
\begin{equation}
 H(r)=\int_0^r \bigl[z(r)q(s)-q(r)z(s)\bigr]s f(s)\,ds.
 \label{eq:corrector-kernel}
\end{equation}
For $0<s\le r$, the kernel has absolute value at most
$\log(r/s)+2$. When $r\le1$, direct integration gives a bound
$C K r^2$; when $r\ge1$, use
$\log(r/s)\le\log r+|\log s|$ and
\[
 \int_0^\infty s(1+|\log s|)(1+t(s)^m)e^{-\alpha t(s)}\,ds<\infty.
\]
This proves \eqref{eq:corrector-linear-bound}. Applying it first to
\eqref{eq:corrector-s} and then to \eqref{eq:corrector-z} proves the
stated growth bounds and uniqueness of the regular radial solutions.

The functions $S,Z$ are $w_0,z_0$ in
\cite[Lemma 9, equations (26), (29)--(31) and Corollary 18]{ManciniMartinazzi}.
The computations in \cite[proof of Proposition 12]{ManciniMartinazzi}
give
\begin{align*}
 \int_{\R^2}(-\Delta S)\,dy&=4\pi,\qquad
 \int_{\R^2}(-\Delta Z)\,dy=12\pi+\frac{2\pi^3}{3},\\
 4\int_{\R^2} e^{2\psi}(S+2\psi S+\psi^2+\psi^3)\,dy
   &=-8\pi-\frac{2\pi^3}{3}.
\end{align*}
The first identity in \eqref{eq:corrector-integrals} follows from
\eqref{eq:corrector-s} and $\int_{\R^2} e^{2\psi}\psi\,dy=-\pi$.
For the second, subtract the bracket in \eqref{eq:corrector-z} from
$B-2Z$ and use the last two displayed identities. All these integrals
converge absolutely by the growth bounds.
\end{proof}

\Needspace{6\baselineskip}
\begin{lemma}
\label{lem:quadratic-concentration-parameter}
Let $u_n$ be centred radial extremals at parameters
$\lambda_n\ge\lambda_0>4\pi$. Suppose
\begin{equation}
 \|u_n\|_2^2\to0,\qquad \int_{\R^2}(e^{4\pi u_n^2}-1)\dd\to C_{\mathrm{cl}}.
 \label{eq:quadratic-concentration-input}
\end{equation}
Put $v_n=\sqrt{4\pi}\,u_n$, $\gamma_n=v_n(0)$ and let $\mu_n$ and
$\mathcal E_n$ be as in \eqref{eq:quadratic-multiplier}. Then
\begin{align}
 \gamma_n^2\|u_n\|_2^2&\longrightarrow1,&
 \frac{\mu_n}{\gamma_n^2}&\longrightarrow C_{\mathrm{cl}},
 \label{eq:quadratic-concentration-scales}\\
 \gamma_n v_n&\longrightarrow4\pi G
     \quad\text{strongly in }L^q(\R^2),\quad 1\le q<\infty,
 \label{eq:quadratic-outer-profile}\\
 \mathcal E_n&=1+
 \left(1+\frac{4\pi}{C_{\mathrm{cl}}}\right)\gamma_n^{-4}
                  +o(\gamma_n^{-4}),
 \label{eq:quadratic-sharp-energy}\\
 \lambda_n&\longrightarrow4\pi+C_{\mathrm{cl}}.
 \label{eq:quadratic-concentration-parameter}
\end{align}
In particular, no such family can have $\lambda_n\to+\infty$.
\end{lemma}

\begin{proof}
Write $S_n=S(\lambda_n,2)$ and $m_n=\int_{\R^2}v_n^2\,dx=4\pi \|u_n\|_2^2$.
The barriers give $C_{\mathrm{cl}}\le S_n\le \int_{\R^2}(e^{4\pi u_n^2}-1)\dd$, so
\begin{equation}
 S_n\to C_{\mathrm{cl}},\quad \lambda_n\|u_n\|_2^2\to0,\quad
 \mu_n\to\infty,\quad
 \frac{\lambda_n}{\mu_n}
   =\frac{\lambda_n\|u_n\|_2^2}{S_n}\to0.
 \label{eq:quadratic-penalty-ratio}
\end{equation}
If $\gamma_n$ were bounded along a subsequence, the estimate
$e^{v_n^2}-1\le C v_n^2$ would contradict \eqref{eq:quadratic-concentration-input}.
Thus $\gamma_n\to\infty$. The pointwise inequality
\[
 \frac{v_n^2e^{v_n^2}}{\gamma_n^2}
 \le\frac{e^{v_n^2}-1}{1-e^{-\gamma_n^2}}
\]
and $\lambda_n>0$ give
\begin{equation}
 \limsup\frac{\mu_n}{\gamma_n^2}\le C_{\mathrm{cl}}.
 \label{eq:quadratic-multiplier-upper}
\end{equation}

\emph{The inner limit and truncation.}
Set
\[
 r_n^2=\frac{\mu_n}{\pi\gamma_n^2}e^{-\gamma_n^2},\quad
 a_n=1+\frac{\lambda_n}{\mu_n},\quad
 \eta_n(y)=\gamma_n\bigl(v_n(r_n y)-\gamma_n\bigr).
\]
Then $r_n^2\le C e^{-\gamma_n^2}$ and $a_n\to1$. Equation~\eqref{eq:quadratic-euler} becomes
\begin{equation}
 -\Delta\eta_n
 =4\left(1+\frac{\eta_n}{\gamma_n^2}\right)
     e^{2\eta_n+\eta_n^2/\gamma_n^2}
 -a_n\gamma_n^2r_n^2
       \left(1+\frac{\eta_n}{\gamma_n^2}\right).
 \label{eq:quadratic-inner-equation}
\end{equation}
Radial monotonicity gives $\eta_n\le0$ and
$0\le1+\eta_n/\gamma_n^2\le1$. In particular
$-\Delta\eta_n\le4$, whence $\eta_n(y)\ge-|y|^2$ by radial integration.
Local elliptic estimates and uniqueness of the regular radial initial-value
problem give
\begin{equation}
 \eta_n\to\psi=-\log(1+|y|^2)
                     \quad\text{in }C^1_{\mathrm{loc}}(\R^2).
 \label{eq:quadratic-inner-limit}
\end{equation}

For $A>1$, let $v_{n,A}=\min\{v_n,\gamma_n/A\}$.
Testing $-\Delta v_n+a_n v_n=(4\pi/\mu_n)v_ne^{v_n^2}$ by
$(v_n-\gamma_n/A)_+$, and integrating first over $B_{Rr_n}$,
gives
\[
 \liminf \|\nabla(v_n-\gamma_n/A)_+\|_2^2
 \ge4\left(1-\frac1A\right)\int_{B_R}e^{2\psi}\,dy.
\]
The omitted mass term is bounded by $a_n m_n\to0$.
Let $R\to\infty$ and use $\|\nabla v_n\|_2^2+m_n=4\pi$. Then
\begin{equation}
 \limsup \|\nabla v_{n,A}\|_2^2\le\frac{4\pi}{A},\qquad
 \int_{\R^2}v_{n,A}^2\,dx\le m_n\to0.
 \label{eq:quadratic-truncation-energy}
\end{equation}
For $1<q<A$, Ruf's inequality and spatial dilation therefore give
\begin{equation}
 \int_{\R^2}\bigl(e^{qv_{n,A}^2}-1\bigr)\,dx
                    \le C_{A,q}\int_{\R^2} v_{n,A}^2\,dx.
 \label{eq:quadratic-truncation-exponential}
\end{equation}
Indeed, choose $\alpha<1$ such that $q\|\nabla v_{n,A}\|_2^2/(4\pi)\le\alpha$.
Dilate $\sqrt{q/(4\pi)}\,v_{n,A}$ to mass $1-\alpha$, apply Ruf's
full-norm bound, and undo the dilation.

Splitting $\int_{\R^2}(e^{4\pi u_n^2}-1)\dd$ at the level $\gamma_n/A$ gives
\[
 \int_{\R^2}(e^{4\pi u_n^2}-1)\dd\le o(1)+\frac{A^2}{\gamma_n^2}\int_{\R^2} v_n^2e^{v_n^2}\,dx
 =o(1)+\frac{A^2}{\gamma_n^2}
                  \bigl(\mu_n+\lambda_n\|u_n\|_2^2\bigr).
\]
Let $A\downarrow1$ and combine with \eqref{eq:quadratic-multiplier-upper}.
This proves $\mu_n/\gamma_n^2\to C_{\mathrm{cl}}$.

\emph{The outer profile.}
Define the nonnegative source
\[
 F_n=\frac{4\pi\gamma_n}{\mu_n}v_n(e^{v_n^2}-1),\qquad
 b_n=1+\frac{\lambda_n-4\pi}{\mu_n}\to1.
\]
Then $(-\Delta+b_n)(\gamma_n v_n)=F_n$. On
$\{v_n>\gamma_n/A\}$,
\[
 \int_{\{v_n>\gamma_n/A\}}F_n\,dx
 \le\frac{4\pi A}{\mu_n}\int_{\R^2} v_n^2e^{v_n^2}\,dx
                      =4\pi A+o(1).
\]
On the complementary set, choose $1<p<\min\{A,2\}$ and
$q=p/(p-1)>2$. By \eqref{eq:quadratic-truncation-exponential},
\[
 \int_{\{v_n\le\gamma_n/A\}}F_n\,dx
 \le\frac{4\pi\gamma_n}{\mu_n}\|v_n\|_q
                  \|e^{v_{n,A}^2}-1\|_p=o(1).
\]
Here the $L^q$ norms of $v_n$ are bounded by their full $H^1$ norms,
and $(e^s-1)^p\le e^{ps}-1$ for $s\ge0$.
The inner limit gives, for each fixed $R$,
\[
 \int_{B_{Rr_n}}F_n\,dx\longrightarrow4\int_{B_R}e^{2\psi}\,dy.
\]
Let $R\to\infty$ and $A\downarrow1$. It follows that
$\int_{\R^2}F_n\,dx\to4\pi$ and $F_n\,dx\rightharpoonup4\pi\delta_0$,
with all its mass concentrating at the origin.

Let $G_b(x)=G(\sqrt b\,x)$, the unit-flux Green kernel for
$-\Delta+b$. Since $G_b\in L^q(\R^2)$ for every finite $q\ge1$,
the resolvent representation
\[
 \gamma_n v_n=G_{b_n}*F_n
\]
and continuity of translations in $L^q$ prove
\eqref{eq:quadratic-outer-profile}. More explicitly, integrate
$\|G_{b_n}(\cdot-y)-G\|_q$ against $F_n(y)$, use uniform continuity
for small $|y|$, and the vanishing source mass outside any fixed ball.
The resolvent identity follows from uniqueness of the $L^2$ solution of
the homogeneous equation. Since $\int_{\R^2}G^2\,dx=1/(4\pi)$, taking $q=2$
gives $\gamma_n^2m_n\to4\pi$, and hence the first limit in
\eqref{eq:quadratic-concentration-scales}.

\emph{The fourth-order inner expansion.}
Fix $\delta=3/4$ and write $\varepsilon_n=\gamma_n^{-2}$.
On $0\le t\le\delta\gamma_n^2$, we claim
\begin{equation}
 \eta_n=\psi+\varepsilon_nS+\varepsilon_n^2 Z
                      +O\bigl(\varepsilon_n^3(1+t)\bigr).
 \label{eq:quadratic-fourth-inner}
\end{equation}
The error is uniform on this interval. Write
$\eta_n=\psi+\varepsilon_nS+\varepsilon_n^2(Z+H_n)$.
On any interval where $|H_n|\le1+t$, Taylor's theorem in
\eqref{eq:quadratic-inner-equation}, after subtracting
\eqref{eq:corrector-s} and \eqref{eq:corrector-z}, gives
\begin{equation}
 \mathcal{L}H_n=f_n,\qquad
 |f_n|\le C\varepsilon_n(1+t^6)e^{-(2-\delta)t}
                       +C\gamma_n^6r_n^2.
 \label{eq:quadratic-bootstrap-source}
\end{equation}
The two nonlinear coefficients subtracted here are
$2S+\psi^2+\psi$ and
$2Z+S+2S^2+4\psi S+2\psi^2S+\psi^3+\psi^4/2$.
At order $\varepsilon_n^2$, the term $8e^{2\psi}H_n$ is moved to the left-hand side and included in $\mathcal{L}H_n$. Every remaining term containing $H_n$ carries an additional factor $\varepsilon_n$ or a difference $e^{X_n}-1$. The bootstrap bound $|H_n|\le1+t$ and the inequality
\[
 |e^X-1|\le |X|e^{\max\{X,0\}}
\]
therefore bound these terms by the first expression in \eqref{eq:quadratic-bootstrap-source}.
To check the Taylor remainder on this growing interval, put
$X_n=2(\eta_n-\psi)+\varepsilon_n\eta_n^2$. Under the assumed
bound on $H_n$,
\[
 X_n=\varepsilon_n(2S+\psi^2)
       +\varepsilon_n^2\bigl[2(Z+H_n)+2\psi S\bigr]
       +O\bigl(\varepsilon_n^3(1+t)^2\bigr).
\]
Moreover $|X_n|\le C\varepsilon_n(1+t)^2$ and
$X_n\le\delta t+C$, since $\varepsilon_nt^2\le\delta t$
and $\varepsilon_nt\le\delta$. Taylor's integral remainder for
$e^{X_n}$ is therefore bounded by
$C\varepsilon_n^3(1+t^6)e^{\delta t}$. Multiplication by
$e^{2\psi}$ and division by $\varepsilon_n^2$ give the first
term in \eqref{eq:quadratic-bootstrap-source}.
The linear term in \eqref{eq:quadratic-inner-equation}, divided by
$\varepsilon_n^2$, is bounded by $C\gamma_n^6r_n^2$.
Moreover,
\[
 \gamma_n^6r_n^2
 =o\bigl(\varepsilon_n e^{-(2-\delta)\delta\gamma_n^2}\bigr),
\]
since $1-(2-\delta)\delta=(1-\delta)^2>0$.
It can therefore be absorbed into the first bound in
\eqref{eq:quadratic-bootstrap-source}. \reftext{Lemma}{lem:quadratic-correctors},
with $\alpha=2-\delta>1$, gives
$|H_n|\le C\varepsilon_n(1+t)$. Since $H_n(0)=H_n'(0)=0$,
the improved bound and a continuation argument extend the estimate
to the whole interval. This proves \eqref{eq:quadratic-fourth-inner}.

Put $\rho_n=r_n\sqrt{e^{\delta\gamma_n^2}-1}$. In inner coordinates,
\eqref{eq:quadratic-fourth-inner} gives
\begin{align*}
 4(1+\varepsilon_n\eta_n)^2e^{2\eta_n+\varepsilon_n\eta_n^2}
 &=4e^{2\psi}(1+\varepsilon_n A+\varepsilon_n^2 B)\\
 &\quad+O\bigl(\varepsilon_n^3(1+t^6)e^{-(2-\delta)t}\bigr).
\end{align*}
The remainder is integrable uniformly, since radial area measure is
$\pi e^t\,dt$ and $2-\delta>1$. The omitted tails of the corrector
integrals are exponentially small. Consequently
\begin{equation}
 \frac{4\pi}{\mu_n}\int_{B_{\rho_n}}v_n^2e^{v_n^2}\,dx
                  =4\pi+4\pi\gamma_n^{-4}+O(\gamma_n^{-6}).
 \label{eq:quadratic-core-energy}
\end{equation}

\emph{The whole-plane remainder and the parameter.}
At $\rho_n$, expansion \eqref{eq:quadratic-fourth-inner} gives
$v_n(\rho_n)=\gamma_n/4+O(\gamma_n^{-1})$. Thus
$v_n\le\gamma_n/3$ outside $B_{\rho_n}$ for large $n$.
The outer profile, H\"older's inequality and
\eqref{eq:quadratic-truncation-exponential} with $A=3$, $q=2$ yield
\begin{align*}
 \frac{4\pi}{\mu_n}\int_{\R^2\setminus B_{\rho_n}}
                         v_n^2(e^{v_n^2}-1)\,dx
 &\le\frac{4\pi}{\mu_n\gamma_n^2}
        \|\gamma_n v_n\|_4^2\|e^{v_{n,3}^2}-1\|_2\\
 &=O(\gamma_n^{-5})=o(\gamma_n^{-4}).
\end{align*}
By \eqref{eq:quadratic-concentration-scales},
$\mu_n\simeq\gamma_n^2$ and $m_n\simeq\gamma_n^{-2}$.
The inequality $(e^s-1)^2\le e^{2s}-1$ and
\eqref{eq:quadratic-truncation-exponential} therefore give
$\|e^{v_{n,3}^2}-1\|_2=O(m_n^{1/2})=O(\gamma_n^{-1})$.
Also $\rho_n^2\le C e^{-(1-\delta)\gamma_n^2}$, so the plain
$L^2$ contribution inside $B_{\rho_n}$ is exponentially small after
multiplication by $4\pi/\mu_n$. Adding the plain mass outside the
inner ball to \eqref{eq:quadratic-core-energy}, we obtain
\begin{align*}
 4\pi\mathcal E_n
 &=4\pi+4\pi\gamma_n^{-4}
           +\frac{4\pi}{\mu_n}\int_{\R^2} v_n^2\,dx+o(\gamma_n^{-4})\\
 &=4\pi+\left(4\pi+\frac{16\pi^2}{C_{\mathrm{cl}}}\right)\gamma_n^{-4}
                                                   +o(\gamma_n^{-4}).
\end{align*}
This proves \eqref{eq:quadratic-sharp-energy}. Finally use the exact
identity in \eqref{eq:quadratic-energy-identity}:
\[
 \lambda_n S_n=\mu_n^2(\mathcal E_n-1)
                   \longrightarrow C_{\mathrm{cl}}(C_{\mathrm{cl}}+4\pi).
\]
Since $S_n\to C_{\mathrm{cl}}>0$, \eqref{eq:quadratic-concentration-parameter}
follows. The argument applies without assuming that $\lambda_n$ is
bounded.
\end{proof}

\section{Attainment thresholds and endpoint compactness}\label{sec:endpoint-proofs}

We begin with the quadratic penalty.

\begin{proof}[Proof of Theorem~\ref{thm:quadratic}]
Write $Q_-(u)$ and $Q_+(u)$ for the two quotients in
\eqref{eq:quadratic-quotients}. The small-gradient bound
\eqref{eq:small-r} with $r=2$ and Ruf's inequality away from zero
gradient energy give $D_*<\infty$.
Lemma~\ref{lem:quadratic-explicit-gap} gives
$D_*>\max\{4\pi,8\pi^2 B_2\}$, and
Proposition~\ref{prop:quadratic-green-test} gives
$\lambda_+>4\pi+C_{\mathrm{cl}}$.
For a normalized $u$, put $d=4\pi-\lambda$. Then
\[
 J_{\lambda,2}(u)-d=\|\nabla u\|_2^2\bigl(Q_-(u)-d\bigr),\qquad
 J_{\lambda,2}(u)-C_{\mathrm{cl}}
    =\|u\|_2^2\bigl(Q_+(u)-\lambda\bigr).
\]
Spatial dilation to full norm increases both quotients whenever the
norm is strictly less than one. Together with the lower bounds $S(\lambda,2)\ge\max\{d,C_{\mathrm{cl}}\}$,
these identities give $S(\lambda,2)>d$ exactly for $\lambda>\lambda_-$,
and $S(\lambda,2)>C_{\mathrm{cl}}$ exactly for $\lambda<\lambda_+$.
Thus Lemma~\ref{lem:quadratic-compactness}, applied with the corrected
barrier $C_{\mathrm{cl}}$, gives interior extremals that we can follow
to the two endpoints.

Take a normalized centred radial maximizing sequence for $Q_-$.
Such sequences exist by rearrangement and spatial dilation to full norm.
The first identity gives $J_{\lambda_-,2}(u_n)\to D_*$.
Since $\lambda_-<0$ and $\|u_n\|_2^2\le1$,
\[
 \liminf_{n\to\infty}\int_{\R^2}(e^{4\pi u_n^2}-1)\dd
 \ge D_*+\lambda_-=4\pi.
\]
A weak zero limit would give $g_n=\|\nabla u_n\|_2^2\to0$
by Lemma~\ref{lem:weakzero}. Then \eqref{eq:small-four} and
$\|u_n\|_4^4\le B_2 g_n(1-g_n)$ give
$\limsup_n Q_-(u_n)\le8\pi^2B_2<D_*$.
Hence Lemma~\ref{lem:quadratic-compactness} gives strong subsequential
convergence to an extremal at $\lambda_-$, which attains $Q_-=D_*$.
Every lower endpoint extremal has this quotient value. For extremals
at $\lambda_n>\lambda_-$ tending to $\lambda_-$ we have
$4\pi-\lambda_n\le Q_-(u_n)\le D_*$ by the first identity. Thus the
same compactness applies to both families of centred radial extremals.

Choose centred radial extremals $u_n$ at interior parameters tending
to $\lambda_+$, allowing initially $\lambda_+=+\infty$. In the infinite
case, $\|u_n\|_2^2\le S(0,2)/\lambda_n\to0$ for large $n$, so every weak
limit is zero. In the finite case, suppose that a subsequence has
weak limit zero. If its gradient energies tend to $g<1$, remainder
continuity gives
$J_{\lambda_n,2}(u_n)\to(4\pi-\lambda_+)(1-g)\le0$,
contrary to $S(\lambda_n,2)\ge C_{\mathrm{cl}}$.
Thus a weakly null subsequence has $\|u_n\|_2^2\to0$ in either case.
The exponential integral is bounded below by $C_{\mathrm{cl}}$ since
$\lambda_n>0$ eventually. The weak-zero alternative in the proof of
Lemma~\ref{lem:quadratic-compactness} therefore forces concentration
and convergence of this integral to $C_{\mathrm{cl}}$.
Lemma~\ref{lem:quadratic-concentration-parameter} would then give
$\lambda_n\to4\pi+C_{\mathrm{cl}}$, contrary to the strict Green test.
It follows that $\lambda_+$ is finite and the weak limits are nonzero.
Lemma~\ref{lem:quadratic-compactness} gives strong subsequential
convergence to an upper endpoint extremal. The same exclusion of a
zero weak limit at the fixed parameter $\lambda_+$ proves compactness
of its centred radial extremals. Since $S(\lambda_+,2)=C_{\mathrm{cl}}$,
each such extremal attains $Q_+=\lambda_+$.

Rearrangement preserves extremality and mass by Lemma~\ref{lem:normalization}.
Endpoint compactness therefore gives $m_-,m_+\in(0,1)$ as an attained
minimum and maximum, respectively. Let $u_h^-$ and $u_h^+$ be extremals
at $\lambda_-+h$ and $\lambda_+-h$. Comparing them with endpoint
extremals of the selected masses gives
\begin{align*}
 D_*-hm_- &\le S(\lambda_-+h,2)\le D_*-h\|u_h^-\|_2^2,\\
 C_{\mathrm{cl}}+hm_+ &\le S(\lambda_+-h,2)
                  \le C_{\mathrm{cl}}+h\|u_h^+\|_2^2.
\end{align*}
Every limiting mass belongs to an endpoint extremal by compactness of the
centred representatives. The first comparison forces $\|u_h^-\|_2^2\to m_-$
and the second forces $\|u_h^+\|_2^2\to m_+$. Substitution gives
\eqref{eq:quadratic-lower-slope} and \eqref{eq:quadratic-upper-slope}.
The exterior slopes follow from the formulas for $S$ outside the interval.
\end{proof}

\begin{proof}[Proof of Theorem~\ref{thm:endpoint}]
Let $\lambda_n<\lamstar$ tend to $\lamstar$. Take centred radial
extremals $u_n$ at $\lambda_n$ using \cite[Theorem 1.1]{Chen} and
Lemma~\ref{lem:normalization}. They have full norm, and
their values tend to $S(\lamstar,p)=4\pi$ by \eqref{eq:lipschitz}.
If a weak limit were zero, then
$\|u_n\|_p^p\to0$ and $\int_{\R^2}(e^{4\pi u_n^2}-1)\dd\to4\pi$.
\reftext{Lemma}{lem:weakzero} gives $g_n=\|\nabla u_n\|_2^2\to0$.
For $2<p<4$, estimate \eqref{eq:small-r} with $r=p$ gives
\[
 J_{\lambda_n,p}(u_n)-4\pi
 \le -4\pi g_n+
    \bigl(Cg_n^{(4-p)/2}-\lambda_n\bigr)\|u_n\|_p^p<0
\]
for large $n$, since $\lambda_n\to\lamstar>0$. This contradicts
\eqref{eq:attained-comparison}.
For $p=4$, \reftext{Proposition}{prop:quartic-gap}
gives $\lambda_n>\lambda_{\mathrm{loc}}+\varepsilon$ eventually for
some $\varepsilon>0$. Using \eqref{eq:small-four} and
$\|u_n\|_4^4\le B_2g_n(1-g_n)$, one obtains
\[
 J_{\lambda_n,4}(u_n)-4\pi
 \le\left[-4\pi+(8\pi^2-\lambda_n)_+B_2+O(g_n)\right]g_n<0,
\]
again a contradiction.

The weak limit is nonzero. \reftext{Lemma}{lem:nonzero-compactness} gives
strong $H^1$ convergence along a subsequence, and the limit attains
$S(\lamstar,p)=4\pi$. If instead $u_n$ are centred radial endpoint
extremals, the same weak-zero argument would contradict their exact
value $4\pi$. Hence the endpoint extremal set is compact as well.
The quotient formula \eqref{eq:threshold-quotient} follows from
$J_{\lamstar,p}(u)\le4\pi$ for every $u\in\A$, with equality at a
nonzero endpoint extremal.

By \reftext{Lemma}{lem:normalization}, the maximum in
\eqref{eq:sigma} may be taken over the compact centred radial endpoint
set. It is attained and is strictly positive.

Choose an endpoint extremal $v$ with $\|v\|_p^p=\sigma_p$. If $u_h$ is
any extremal at $\lamstar-h$, then
\begin{equation}
 \sigma_p\le\frac{S(\lamstar-h,p)-4\pi}{h}\le \|u_h\|_p^p.
 \label{eq:derivative-squeeze}
\end{equation}
The first inequality follows by testing $v$ at $\lamstar-h$,
and the second follows from
$J_{\lamstar,p}(u_h)\le4\pi$.
For every sequence $h_n\downarrow0$, the centred rearrangements of
$u_{h_n}$ have a subsequence converging strongly to an endpoint
extremal $w$. Their power masses converge to $\|w\|_p^p\le\sigma_p$.
Combining with \eqref{eq:derivative-squeeze} proves both the limit
of the difference quotient and $\|u_h\|_p^p\to\sigma_p$ along every
family. It also identifies the selected endpoint subset.
\end{proof}

\begin{remark}
\label{rem:endpoint-sequences}
For a fixed unit-mass profile, the family \eqref{eq:test-family} has
$\|\nabla U_{t,v}\|_2^2\to0$, $\|U_{t,v}\|_2^2\to1$ and weak limit zero as $t\downarrow0$.
It maximizes the endpoint functional asymptotically: the value tends
to $4\pi$ for $2<p\le4$ at $\lamstar$, and to $D_*$ for $p=2$ at
$\lambda_-$. The strict lower gap excludes it from maximizing the
quotient in \eqref{eq:quadratic-quotients}.
At $\lambda_+$, there also exist normalized concentrating sequences
with $\int_{\R^2}(e^{4\pi u_n^2}-1)\dd\to C_{\mathrm{cl}}$ and $\|u_n\|_2^2\to0$, hence
$J_{\lambda_+,2}(u_n)\to C_{\mathrm{cl}}$; the Green test family
\eqref{eq:green-test} provides such a sequence.
Thus centred radial endpoint extremals are compact, although the
endpoint functionals also admit noncompact maximizing sequences.
\end{remark}

\section{Asymptotics and profile selection}\label{sec:strong-proofs}

We first establish compactness for the limiting problem at unit mass.

\Needspace{6\baselineskip}
\begin{lemma}\label{prop:limit-attainment}
For $p>4$, $0<\Cp<\infty$. The supremum in
\eqref{eq:limit-problem} is achieved by a centred radial function,
and every centred radial maximizing sequence with unit $L^2$ mass
has a strongly convergent subsequence in $H^1$.
\end{lemma}

\begin{proof}
For $\|v\|_2=1$, interpolation between $L^2$ and $L^p$ yields
\begin{equation}
 \|v\|_4^4\le \|v\|_p^{2p/(p-2)}.
 \label{eq:mass-interpolation}
\end{equation}
Therefore
\[
 \Hfun(v)\le8\pi^2\|v\|_p^{2p/(p-2)}-\|v\|_p^p,
\]
whose right-hand side is bounded above because $2/(p-2)<1$. For positivity, choose
$v_0$ from \eqref{eq:rational-test} and set $v_s(x)=s v_0(sx)$.
Then $\|v_s\|_2=1$ and
\[
 \Hfun(v_s)=s^2\frac{24\pi}{35}-s^{p-2}\|v_0\|_p^p>0
\]
for small positive $s$, since $p>4$.

Schwarz rearrangement does not decrease $\Hfun$, so we may choose a
centred radial maximizing sequence $w_n$ with $\|w_n\|_2=1$.
Its positive values imply
\[
 4\pi \|\nabla w_n\|_2^2+\|w_n\|_p^p\le8\pi^2\|w_n\|_p^{2p/(p-2)}
\]
for all sufficiently large $n$. Hence both $\|\nabla w_n\|_2^2$ and $\|w_n\|_p^p$
are bounded. Take a weak $H^1$ limit $w$. The radial tail argument
following \eqref{eq:radial-bound} gives strong $L^4$ and $L^p$
convergence. Weak lower semicontinuity then gives
\begin{equation}
 \Hfun(w)\ge\limsup_n\Hfun(w_n)=\Cp>0,
 \qquad 0<\|w\|_2^2=m\le1.
 \label{eq:limit-lsc}
\end{equation}
If $m<1$, set $z(x)=w(\sqrt m\,x)$. Then $\|z\|_2=1$,
$\|\nabla z\|_2^2=\|\nabla w\|_2^2$, and
\[
 \Hfun(z)=m^{-1}\Hfun(w)+4\pi(m^{-1}-1)\|\nabla w\|_2^2
 >\Hfun(w)\ge\Cp,
\]
which contradicts the definition of $\Cp$. Thus $m=1$.
Equation~\eqref{eq:limit-lsc} is now an equality. Convergence of the
quartic and power moments gives convergence of the gradient norms,
and $\|w_n\|_2=\|w\|_2=1$ gives strong $L^2$ convergence. Hence
$w_n\to w$ strongly in $H^1$. This proves precompactness for every centred radial maximizing
sequence.
\end{proof}

\Needspace{6\baselineskip}
\begin{lemma}\label{prop:identities}
If $v$ maximizes \eqref{eq:limit-problem}, then
\begin{equation}
 \|v\|_p^p=\frac{2\Cp}{p-4},\qquad
 8\pi^2\|v\|_4^4-4\pi \|\nabla v\|_2^2=\frac{p-2}{p-4}\Cp.
 \label{eq:limit-identities}
\end{equation}
A nonnegative maximizer satisfies, in the weak sense,
\begin{equation}
 -4\pi\Delta v+\eta v=16\pi^2v^3-\frac p2v^{p-1},
 \qquad\eta=\Cp+4\pi \|\nabla v\|_2^2>0.
 \label{eq:limit-equation}
\end{equation}
\end{lemma}

\begin{proof}
The mass-preserving family $v_s(x)=s v(sx)$ is admissible in
\eqref{eq:limit-problem}, and
\[
 \Hfun(v_s)=s^2\bigl[8\pi^2\|v\|_4^4-4\pi \|\nabla v\|_2^2\bigr]
 -s^{p-2}\|v\|_p^p.
\]
Its derivative at $s=1$ vanishes. Thus
$2[8\pi^2\|v\|_4^4-4\pi \|\nabla v\|_2^2]=(p-2)\|v\|_p^p$.
Combine this with $\Hfun(v)=\Cp$ to obtain
\eqref{eq:limit-identities}.

The constraint $\|v\|_2=1$ is regular and the functional is continuously
differentiable on $H^1$ for finite $p>4$. The Lagrange multiplier
equation is \eqref{eq:limit-equation} with some real $\eta$.
Testing it by $v$ and using \eqref{eq:limit-identities} gives
$\eta=\Cp+4\pi \|\nabla v\|_2^2>0$.
\end{proof}

Thus all limiting maximizers have the same $L^p$ mass.

These compactness results allow us to identify the rescaled limits of spreading extremals.

\begin{proof}[Proof of Theorem~\ref{thm:strong}]

Take centred radial extremals $u_\lambda$ of full norm and let
$\lambda\to+\infty$. From \eqref{eq:attained-comparison} and the definition of $S(0,p)$,
\[
 \lambda\|u_\lambda\|_p^p
 =\int_{\R^2}(e^{4\pi u_\lambda^2}-1)\dd-S(\lambda,p)
 \le S(0,p)-4\pi.
\]
Thus $\|u_\lambda\|_p^p\to0$ and every weak limit is zero. Since the
exponential integral is at least $S(\lambda,p)>4\pi$ for $\lambda>0$,
\reftext{Lemma}{lem:weakzero} implies
\begin{equation}
 g_\lambda:=\|\nabla u_\lambda\|_2^2\longrightarrow0,\qquad
 \|u_\lambda\|_2^2=1-g_\lambda\longrightarrow1.
 \label{eq:strong-vanishing}
\end{equation}
By \eqref{eq:small-four},
\begin{equation}
 S(\lambda,p)-4\pi
 =-4\pi g_\lambda+8\pi^2\|u_\lambda\|_4^4-\lambda \|u_\lambda\|_p^p
 +O(g_\lambda \|u_\lambda\|_4^4).
 \label{eq:uniform-upper-expansion}
\end{equation}
Interpolation gives
\begin{equation}
 \|u_\lambda\|_4^4\le
 \|u_\lambda\|_2^{2(p-4)/(p-2)}\|u_\lambda\|_p^{2p/(p-2)}
 \le \|u_\lambda\|_p^{2p/(p-2)}.
 \label{eq:original-interpolation}
\end{equation}
The strict positivity of the left side of
\eqref{eq:uniform-upper-expansion} yields
\[
 \lambda \|u_\lambda\|_p^p<(8\pi^2+Cg_\lambda)\|u_\lambda\|_4^4,
 \qquad
 4\pi g_\lambda<(8\pi^2+Cg_\lambda)\|u_\lambda\|_4^4.
\]
With $\delta=\lambda^{-1/(p-4)}$, these inequalities and
\eqref{eq:original-interpolation} give
\begin{equation}
 \|u_\lambda\|_p^p=O(\delta^{p-2}),\qquad
 \|u_\lambda\|_4^4=O(\delta^2),\qquad
 g_\lambda=O(\delta^2).
 \label{eq:scale-bounds}
\end{equation}
In particular, the remainder in
\eqref{eq:uniform-upper-expansion} is $O(\delta^4)$.

Fix $v\in H^1(\R^2)$ with $\|v\|_2=1$ and use $U_{\delta^2,v}$
from \eqref{eq:test-family}.
It has full $H^1$ norm. Since $\lambda\delta^{p-2}=\delta^2$,
the identities in \eqref{eq:test-moments} and the uniform
remainder estimate give
\begin{align}
 \frac{J_{\lambda,p}(U_{\delta^2,v})-4\pi}{\delta^2}
 &=\frac{-4\pi \|\nabla v\|_2^2}{1+\delta^2\|\nabla v\|_2^2}
 +\frac{8\pi^2\|v\|_4^4}{(1+\delta^2\|\nabla v\|_2^2)^2}
 \nonumber\\
 &\quad-\frac{\|v\|_p^p}{(1+\delta^2\|\nabla v\|_2^2)^{p/2}}+O(\delta^2)
 \nonumber\\
 &\longrightarrow\Hfun(v).
 \label{eq:strong-lower}
\end{align}
Taking the supremum over fixed profiles gives
\begin{equation}
 \liminf_{\lambda\to\infty}
 \delta^{-2}[S(\lambda,p)-4\pi]\ge\Cp.
 \label{eq:liminf}
\end{equation}

Define $v_\lambda$ as in \eqref{eq:profile-scaling}. Then
\begin{equation}
 \begin{aligned}
 \|v_\lambda\|_2^2&=1-g_\lambda,&
 \|\nabla v_\lambda\|_2^2&=\delta^{-2}g_\lambda,\\
 \|v_\lambda\|_4^4&=\delta^{-2}\|u_\lambda\|_4^4,&
 \|v_\lambda\|_p^p&=\delta^{-(p-2)}\|u_\lambda\|_p^p.
 \end{aligned}
 \label{eq:rescaled-moments}
\end{equation}
All moments and the gradient energy in this display are bounded by
\eqref{eq:scale-bounds}. Equation~\eqref{eq:uniform-upper-expansion} becomes
\begin{equation}
 \delta^{-2}[S(\lambda,p)-4\pi]=\Hfun(v_\lambda)+O(\delta^2).
 \label{eq:rescaled-value}
\end{equation}
Normalize the rescaled mass by
\begin{equation}
 w_\lambda=\frac{v_\lambda}{\|v_\lambda\|_2},\qquad \|w_\lambda\|_2=1.
 \label{eq:mass-renormalization}
\end{equation}
The bounded moments in \eqref{eq:rescaled-moments} and $g_\lambda=O(\delta^2)$ give
\begin{equation}
 \Hfun(w_\lambda)-\Hfun(v_\lambda)=O(\delta^2).
 \label{eq:normalization-error}
\end{equation}
Since $\Hfun(w_\lambda)\le\Cp$, \eqref{eq:rescaled-value} and
\eqref{eq:normalization-error} bound the scaled value above by
$\Cp+O(\delta^2)$. Testing a fixed maximizer of \eqref{eq:limit-problem}
in \eqref{eq:strong-lower} gives the lower bound $\Cp-O(\delta^2)$.
Hence
\begin{equation}
 S(\lambda,p)=4\pi+\Cp\delta^2+O(\delta^4).
 \label{eq:first-order-strong}
\end{equation}
The coefficient of $\delta^4$ will be determined below.

Equations~\eqref{eq:rescaled-value} and
\eqref{eq:normalization-error} now imply
$\Hfun(w_\lambda)\to\Cp$. Thus $w_\lambda$ is a centred radial
maximizing family for the unit-mass limiting problem.
\reftext{Lemma}{prop:limit-attainment} gives strong $H^1$
subsequential convergence to a maximizer $v$.
Moreover,
\[
 \|w_\lambda-v_\lambda\|_{H^1}
 =\bigl|(1-g_\lambda)^{-1/2}-1\bigr|
 \|v_\lambda\|_{H^1}=O(\delta^2),
\]
so $v_\lambda$ has the same strong limits. By
\eqref{eq:limit-identities}, each such limit has
$\|v\|_p^p=2\Cp/(p-4)$ and
$8\pi^2\|v\|_4^4-4\pi \|\nabla v\|_2^2=(p-2)\Cp/(p-4)$.
The subsequence principle therefore gives
\eqref{eq:strong-mass} and \eqref{eq:strong-exp} along the whole
centred radial family.

Schwarz rearrangement preserves the gradient energy and nonlinear moments of an extremal by Lemma~\ref{lem:normalization}. Hence the integral limits hold for arbitrary extremal families, while the profile convergence is stated for centred radial representatives. To determine the second coefficient and the selected profiles, we combine compactness of $\mathcal G_p$ with a uniform expansion through the sixth exponential moment.

To compute the normalization correction, define, for $\|v\|_2=1$,
\begin{equation}
 \begin{aligned}
 \mathcal Q_p(v)&=4\pi\|\nabla v\|_2^4
                 -16\pi^2\|\nabla v\|_2^2\|v\|_4^4\\
 &\quad+\frac p2\|\nabla v\|_2^2\|v\|_p^p
       +\frac{32\pi^3}{3}\|v\|_6^6.
 \end{aligned}
 \label{eq:second-functional}
\end{equation}
This functional is continuous on $H^1$. The identities in
\eqref{eq:limit-identities} give
\begin{equation}
 \mathcal Q_p(v)=\frac{32\pi^3}{3}\|v\|_6^6-4\pi \|\nabla v\|_2^4-\Cp \|\nabla v\|_2^2,
 \qquad v\in\mathcal G_p.
 \label{eq:second-coefficient-reduced}
\end{equation}
Consequently, $D_p=\max_{v\in\mathcal G_p}\mathcal Q_p(v)$.
For $\|v\|_2=1$ and $\|\nabla v\|_2^2\le K$, use the test $U_{\delta^2,v}$
defined in \eqref{eq:test-family}.
The moment identities \eqref{eq:test-moments} and the remainder
estimate \eqref{eq:small-six} give, uniformly on this set as
$\delta\downarrow0$,
\begin{align}
 J_{\lambda,p}(U_{\delta^2,v})
 &=\frac{4\pi}{1+\delta^2 \|\nabla v\|_2^2}
   +\frac{8\pi^2\delta^2 \|v\|_4^4}{(1+\delta^2 \|\nabla v\|_2^2)^2}\nonumber\\
 &\quad+\frac{32\pi^3\delta^4\|v\|_6^6}{3(1+\delta^2 \|\nabla v\|_2^2)^3}
   -\frac{\delta^2 \|v\|_p^p}{(1+\delta^2 \|\nabla v\|_2^2)^{p/2}}+O_K(\delta^6)\nonumber\\
 &=4\pi+\delta^2\Hfun(v)+\delta^4\mathcal Q_p(v)+O_K(\delta^6).
 \label{eq:uniform-second-expansion}
\end{align}
Here $\lambda=\delta^{-(p-4)}$. The constant may depend on the fixed $p$.
Indeed, $\|\nabla U_{\delta^2,v}\|_2^2\le K\delta^2$ and
$\|U_{\delta^2,v}\|_6^6=\delta^4(1+\delta^2 \|\nabla v\|_2^2)^{-3}\|v\|_6^6$, while all fixed moments
of $v$ are bounded by its $H^1$ norm. Thus the exponential remainder
in \eqref{eq:small-six} is $O_K(\delta^6)$. Taylor expansion of the four
denominators has the same uniform error. This proves
\eqref{eq:uniform-second-expansion} without a pointwise bound on $v$.

For an extremal $u_\lambda$, the mass normalization gives
\begin{equation}
 \begin{aligned}
 \|\nabla w_\lambda\|_2^2&=\frac{g_\lambda}{\delta^2(1-g_\lambda)},&
 1+\delta^2\|\nabla w_\lambda\|_2^2&=\frac1{1-g_\lambda},\\
 u_\lambda&=U_{\delta^2,w_\lambda}.
 \end{aligned}
 \label{eq:second-exact-normalization}
\end{equation}
The energies of $w_\lambda$ are bounded by \eqref{eq:scale-bounds}.
The same expansion and remainder bound therefore apply along the varying extremal profiles:
\begin{equation}
 \frac{S(\lambda,p)-4\pi-\delta^2\Cp}{\delta^4}
  =-\frac{\Cp-\Hfun(w_\lambda)}{\delta^2}
       +\mathcal Q_p(w_\lambda)+O(\delta^2).
 \label{eq:second-value-identity}
\end{equation}
The deficit $\Cp-\Hfun(w_\lambda)$ is nonnegative.

Choose $v\in\mathcal G_p$ with $\mathcal Q_p(v)=D_p$ and test $U_{\delta^2,v}$.
The left-hand side of \eqref{eq:second-value-identity} has lower
limit at least $D_p$. On the other hand, the precompactness proved above
places every strong limit of $w_\lambda$ in $\mathcal G_p$.
Continuity of $\mathcal Q_p$ then gives
$\limsup_{\lambda\to\infty}\mathcal Q_p(w_\lambda)\le D_p$.
Using the nonnegativity of the deficit in \eqref{eq:second-value-identity}
gives the matching upper bound. Hence the left-hand side converges to
$D_p$, which is \eqref{eq:main-asymptotic}.

The same identity now implies
\begin{equation}
 \mathcal Q_p(w_\lambda)\longrightarrow D_p,\qquad
 \Cp-\Hfun(w_\lambda)=o(\delta^2).
 \label{eq:second-deficit}
\end{equation}
Every strong limit is therefore a maximizer $v\in\mathcal G_p$ with
$\mathcal Q_p(v)=D_p$. Precompactness and the subsequence principle also give
\[
 \operatorname{dist}_{H^1}\!\left(v_\lambda,
 \{v\in\mathcal G_p:\mathcal Q_p(v)=D_p\}\right)\longrightarrow0,
\]
since $\|v_\lambda-w_\lambda\|_{H^1}=O(\delta^2)$.
Together with \eqref{eq:second-coefficient-reduced}, this proves
\eqref{eq:second-profile-selection}.
\end{proof}

Rearrangement also preserves $\mathcal H_p$ and $\mathcal Q_p$ for
the normalized rescaling of an extremal. Hence the scalar conclusions
in \eqref{eq:second-deficit} hold for arbitrary extremal families.

\begin{proof}[Proof of Theorem~\ref{thm:p4-transition}]
Fix $\lambda>\lambda_*(4)$ and choose centred radial extremals
$u_\varepsilon$ at $p=4+\varepsilon$. Their existence follows from
\cite[Theorem 1.1]{Chen}, and Lemma~\ref{lem:normalization} gives
$\|u_\varepsilon\|_{H^1}=1$.

We first show that their gradient energy tends to zero. On every bounded
$H^1$ ball, Sobolev embedding and the mean value theorem give
\begin{equation}
 \left|\|w\|_{4+\varepsilon}^{4+\varepsilon}-\|w\|_4^4\right|
 \le C\varepsilon
 \quad\text{for }0<\varepsilon\le\tfrac12.
 \label{eq:p4-uniform-moments}
\end{equation}
Indeed, for $s\ge0$ one has
$|s^{4+\varepsilon}-s^4|\le C\varepsilon(s^3+s^5)$, and the $L^3$
and $L^5$ moments are uniformly bounded on that ball. Taking suprema on
$\A$ gives $S(\lambda,4+\varepsilon)\to S(\lambda,4)=4\pi$.
Moreover,
\[
 J_{\lambda,4}(u_\varepsilon)
 =S(\lambda,4+\varepsilon)
 +\lambda\bigl(\|u_\varepsilon\|_{4+\varepsilon}^{4+\varepsilon}
                 -\|u_\varepsilon\|_4^4\bigr)\longrightarrow4\pi.
\]
Thus these functions form a maximizing family for the quartic problem.
If a subsequence had a nonzero weak limit, Lemma~\ref{lem:nonzero-compactness}
would give a quartic extremal, contradicting nonattainment at
$\lambda>\lambda_*(4)$. Every weak limit is therefore zero. Since
$\lambda>0$ and
$\int_{\R^2}(e^{4\pi u_\varepsilon^2}-1)\dd\ge
S(\lambda,4+\varepsilon)>4\pi$, Lemma~\ref{lem:weakzero} gives
\begin{equation}
 g_\varepsilon=\|\nabla u_\varepsilon\|_2^2\longrightarrow0,
 \qquad\|u_\varepsilon\|_2^2=1-g_\varepsilon\longrightarrow1.
 \label{eq:p4-vanishing}
\end{equation}

For a lower bound, take the fixed optimizer $V$ from
Lemma~\ref{lem:gn-moment}, with unit mass and unit gradient norm.
The moment identities for $U_{t,V}$ and \eqref{eq:small-four} yield
\begin{equation}
 J_{\lambda,4+\varepsilon}(U_{t,V})-4\pi
 =B_2\lambda_{\mathrm{loc}}t
  -\lambda\|V\|_{4+\varepsilon}^{4+\varepsilon}t^{1+\varepsilon/2}
  +O(t^2),
 \label{eq:p4-test-expansion}
\end{equation}
uniformly for small $\varepsilon$ and $t$. The exponential remainder
is $O(t^2)$, and expanding the penalty denominator introduces an
$O(t^{2+\varepsilon/2})$ error with a uniform constant. Choose
\begin{equation}
 t_\varepsilon=
 \left(\frac{2B_2\lambda_{\mathrm{loc}}}
 {(2+\varepsilon)\lambda\|V\|_{4+\varepsilon}^{4+\varepsilon}}
 \right)^{2/\varepsilon}.
 \label{eq:p4-test-scale}
\end{equation}
The base tends to $\lambda_{\mathrm{loc}}/\lambda<1$, so
$t_\varepsilon/\varepsilon\to0$. At this scale, the first two terms in
\eqref{eq:p4-test-expansion} have sum
$\varepsilon B_2\lambda_{\mathrm{loc}}t_\varepsilon/(2+\varepsilon)$.
Consequently,
\begin{equation}
 S(\lambda,4+\varepsilon)-4\pi
 \ge\frac{\varepsilon}{2+\varepsilon}
 B_2\lambda_{\mathrm{loc}}t_\varepsilon\bigl(1+o(1)\bigr),
 \qquad
 \liminf_{\varepsilon\downarrow0}\varepsilon
 \log\bigl(S(\lambda,4+\varepsilon)-4\pi\bigr)
 \ge2\log\frac{\lambda_{\mathrm{loc}}}{\lambda}.
 \label{eq:p4-lower-rate}
\end{equation}

For the upper bound, \eqref{eq:small-four} gives the uniform expansion
\begin{equation}
 \begin{aligned}
 S(\lambda,4+\varepsilon)-4\pi
 &=-4\pi g_\varepsilon+8\pi^2\|u_\varepsilon\|_4^4
   -\lambda\|u_\varepsilon\|_{4+\varepsilon}^{4+\varepsilon}\\
 &\quad+O\bigl(g_\varepsilon\|u_\varepsilon\|_4^4\bigr).
 \end{aligned}
 \label{eq:p4-extremal-expansion}
\end{equation}
Positivity of the left side and the sharp Gagliardo--Nirenberg inequality
give
\begin{equation}
 c g_\varepsilon\le\|u_\varepsilon\|_4^4
 \le B_2g_\varepsilon(1-g_\varepsilon),\qquad
 0<S(\lambda,4+\varepsilon)-4\pi\le Cg_\varepsilon.
 \label{eq:p4-energy-comparison}
\end{equation}
Interpolation between $L^2$ and $L^{4+\varepsilon}$ also gives
\[
 \|u_\varepsilon\|_{4+\varepsilon}^{4+\varepsilon}
 \ge(1-g_\varepsilon)^{-\varepsilon/2}
       \bigl(\|u_\varepsilon\|_4^4\bigr)^{1+\varepsilon/2}.
\]
Insert this into \eqref{eq:p4-extremal-expansion} and use the upper
quartic bound in \eqref{eq:p4-energy-comparison} to obtain
\begin{align*}
 0&<S(\lambda,4+\varepsilon)-4\pi\\
 &\le\left(8\pi^2-\frac{4\pi}{B_2(1-g_\varepsilon)}
                      +Cg_\varepsilon\right)\|u_\varepsilon\|_4^4
 -\lambda(1-g_\varepsilon)^{-\varepsilon/2}
        \bigl(\|u_\varepsilon\|_4^4\bigr)^{1+\varepsilon/2}.
\end{align*}
The coefficient in parentheses tends to $\lambda_{\mathrm{loc}}>0$.
Dividing by the positive quartic moment and using its lower bound in
\eqref{eq:p4-energy-comparison} therefore yields
\[
 g_\varepsilon\le C\|u_\varepsilon\|_4^4
 \le C(1-g_\varepsilon)
       \left(\frac{\lambda_{\mathrm{loc}}+o(1)}{\lambda}
       \right)^{2/\varepsilon}.
\]
It follows that
$\limsup\varepsilon\log g_\varepsilon
\le2\log(\lambda_{\mathrm{loc}}/\lambda)$.
Conversely, \eqref{eq:p4-lower-rate} and
$S(\lambda,4+\varepsilon)-4\pi\le Cg_\varepsilon$ give the matching
lower bound for $\varepsilon\log g_\varepsilon$. These same inequalities
give the matching upper bound for the logarithm of the excess value.
This proves \eqref{eq:p4-value-rate} and the first limit in
\eqref{eq:p4-scale-rate}. Its second limit follows from
$1-g_\varepsilon\to1$.

Finally, the scaling in \eqref{eq:p4-profile-scaling} gives
\begin{equation}
 \|v_\varepsilon\|_2=\|\nabla v_\varepsilon\|_2=1,\qquad
 \|u_\varepsilon\|_r^r
 =(1-g_\varepsilon)g_\varepsilon^{(r-2)/2}
       \|v_\varepsilon\|_r^r.
 \label{eq:p4-profile-moments}
\end{equation}
Divide \eqref{eq:p4-extremal-expansion} by $g_\varepsilon$ and use
these identities. The already proved rate implies
$g_\varepsilon^{\varepsilon/2}\to\lambda_{\mathrm{loc}}/\lambda$.
Since \eqref{eq:p4-uniform-moments} applies uniformly to $v_\varepsilon$,
whose $H^1$ norm is $\sqrt2$, the resulting expression becomes
\[
 0<\frac{S(\lambda,4+\varepsilon)-4\pi}{g_\varepsilon}
 =-4\pi+\frac{4\pi}{B_2}\|v_\varepsilon\|_4^4+o(1).
\]
Together with $\|v_\varepsilon\|_4^4\le B_2$, this forces
$\|v_\varepsilon\|_4^4\to B_2$.
For any sequence $\varepsilon_n\downarrow0$, radial compactness gives
a subsequence converging weakly in $H^1$ and strongly in $L^4$ to $V$.
The sharp inequality and weak lower semicontinuity give
\[
 B_2=\|V\|_4^4\le B_2\|\nabla V\|_2^2\|V\|_2^2\le B_2.
\]
Thus both limiting norms equal one, and the convergence is strong in
$H^1$. The limit is a centred radial optimizer, as claimed.
Rearrangement preserves the gradient energy and all nonlinear moments
of every extremal by Lemma~\ref{lem:normalization}, so the scalar
conclusions also hold for arbitrary extremal families.
\end{proof}

% We recall the profile correspondence underlying \eqref{eq:transfer},
% established in~\cite[Theorem 1.5(ii)]{GuoLiu}.
Finally, we prove the exact anisotropic correspondence \eqref{eq:transfer}.
Convex symmetrization~\cite{AlvinoEtAl} reduces the anisotropic problem
to centred Wulff-radial functions.
If $v(y)=V(|y|)$ is Euclidean radial, set
$u(x)=\theta^{-1/2}V(F^\circ(x))$. Radial integration and
$F(\nabla F^\circ)=1$ almost everywhere yield
\[
 \begin{aligned}
 \int_{\R^2} F(\nabla u)^2\dd&=\|\nabla v\|_2^2,&
 \|u\|_2&=\|v\|_2,\\
 \|u\|_p^p&=\theta^{1-p/2}\|v\|_p^p,&
 \int_{\R^2}(e^{4\kappa_2 u^2}-1)\dd
 &=\theta\int_{\R^2}(e^{4\pi v^2}-1)\,dy.
 \end{aligned}
\]
Convex symmetrization preserves the exponential and power integrals separately and does not increase the anisotropic energy. Schwarz rearrangement gives the reverse comparison through the inverse profile map. These two comparisons prove \eqref{eq:transfer} and transfer extremals in both directions. They require neither smoothness nor a quadratic structure of $F$.

\begin{proof}[Proof of \reftext{Corollary}{cor:anisotropy}]
Apply \eqref{eq:transfer} and its radial mass identities to
the endpoint and asymptotic formulas of
Theorems~\ref{thm:endpoint} and~\ref{thm:strong}.
The ranges away from the threshold transfer from
Chen--Jiang--Lu--Zhu~\cite[Theorem 1.1]{Chen}. For the value
coefficient the factor is
\[
 \theta\,(\theta^{-p/2})^{-2/(p-4)}
 =\theta^{1+p/(p-4)}=\theta^{2(p-2)/(p-4)}.
\]
For the second value coefficient it is
$\theta(\theta^{-p/2})^{-4/(p-4)}
=\theta^{(3p-4)/(p-4)}$.
The power-mass coefficient has the same factor, since
\[
 \theta^{1-p/2}(\theta^{-p/2})^{-(p-2)/(p-4)}
 =\theta^{2(p-2)/(p-4)}.
\]
As in the Euclidean problem, symmetrization of an extremal preserves
its nonlinear moments and has full norm, so the integral statements
hold without imposing symmetry on the original extremal.
For \eqref{eq:anisotropic-p4-rate}, the transferred penalty
$\lambda\theta^{-p/2}$ tends to $\lambda\theta^{-2}>\lambda_*(4)$.
For each sufficiently small $\eta>0$, it lies eventually between
$\lambda\theta^{-2}-\eta$ and $\lambda\theta^{-2}+\eta$.
Monotonicity sandwiches the Euclidean excess value between those
at these two fixed penalties. Apply Theorem~\ref{thm:p4-transition},
let $\eta\downarrow0$, and use $(p-4)\log\theta\to0$ in
\eqref{eq:transfer} to obtain the asserted limit.
\end{proof}

\begin{proof}[Proof of \reftext{Corollary}{cor:quadratic-transfer}]
At $p=2$, \eqref{eq:transfer} reads
$S_F(\lambda,2)=\theta S(\lambda/\theta,2)$.
Its two barriers are $4\kappa_2-\lambda$ and $\theta C_{\mathrm{cl}}$;
the endpoint values, compactness and slopes follow from
Theorem~\ref{thm:quadratic} and the radial correspondence.
\end{proof}

%The exact thresholds, uniqueness of endpoint extremals and the
%structure of the profiles selected by $D_p$ remain open. Further
%questions concern the prefactor in the logarithmic transition, the
%joint limit near $\lambda_*(4)$, the sign of $D_p$ and a strong-penalty
%expansion uniform as $p\downarrow4$.

\medskip
\noindent\textbf{Acknowledgements.} The authors were supported by the National
Natural Science Foundation of China (Nos.12201089, 12371051), the
Chongqing Normal University Foundation (No.BWQB2023006), and the
Open Research Fund of the Key Laboratory of Nonlinear Analysis and
Applications, Ministry of Education, Central China Normal University
(No.NAA2025ORG006).

\end{document}